\documentclass[11pt,reqno]{amsart}
\usepackage[margin=0.8in]{geometry}
\usepackage{amsmath}
\usepackage{amsfonts}
\usepackage{graphicx}
\usepackage{caption}
\usepackage[labelformat=parens,labelfont={}]{subcaption}
\usepackage{amssymb}

\usepackage{esint}
\usepackage{filecontents}
\usepackage{multirow}  % -- tables

\usepackage{
tikz,
pgfplots,
pgfplotstable,
ifthen,
}
\usepgfplotslibrary{patchplots}
\usetikzlibrary{fit}
\usetikzlibrary{calc}
\pgfplotsset{compat=1.15}

\newcommand{\bs}{\boldsymbol}
\newcommand{\abs}[1]{\left| #1 \right|}
\newcommand{\norm}[1]{\left|\!\left| #1 \right|\!\right|} 
\newcommand{\ol}{\overline}
\newcommand{\dif}{\thinspace\mathrm{d}}
\newcommand{\real}{\mathbb{R}}
\newcommand{\complex}{\mathbb{C}}
\newcommand{\Real}[1]{\text{Re}(#1)}
\newcommand{\Imag}[1]{\text{Im}(#1)}
\newcommand{\iu}{\mathfrak{i}} % imaginary unit
\newcommand{\PML}{\textrm{\fontsize{5}{0}\selectfont PML}}
\newcommand{\CPML}{\textrm{\fontsize{5}{0}\selectfont LWC}}
\newcommand{\ABC}{\textrm{\fontsize{5}{0}\selectfont ABC}}
\newcommand{\FE}{\textrm{\fontsize{5}{0}\selectfont FE}}
\newcommand{\weak}{\textrm{weak}}
\newcommand{\HAT}{\textrm{\tiny hat}}

\newcommand{\nref}{\textrm{ref}}
\newcommand{\rint}{\fint}
\newcommand{\oneD}{\textrm{\fontsize{5}{0}\selectfont 1D}}
\newcommand{\threeD}{\textrm{\fontsize{5}{0}\selectfont 3D}}
\newcommand{\uref}{\widehat{u}}  % solution in ref. coord.
\newcommand{\wref}{\widehat{w}}  % solution in ref. coord.
\newcommand{\vs}{\vspace{0.9\baselineskip}}

\begin{document}
\title[]{Absorbing boundary conditions for the Helmholtz equation using Gauss-Legendre quadrature reduced integrations}
\author{Koki Sagiyama}
\address[Koki Sagiyama]{Department of Mathematics at Imperial College London, 738 Huxley Building, South Kensington Campus, London, SW7 2AZ}
\email{k.sagiyama@imperial.ac.uk}
\date{\today}
\begin{abstract}
We introduce a new class of absorbing boundary conditions (ABCs)
for the Helmholtz equation.
The proposed ABCs are obtained by
using $L$ discrete layers and the $Q_N$ Lagrange finite element
in conjunction with
the $N$-point Gauss-Legendre quadrature reduced integration rule
in a specific formulation of perfectly matched layers.
The proposed ABCs are classified by a tuple $(L,N)$, and achieve
reflection error of order $O(R^{2LN})$ for some $R<1$.
The new ABCs generalise the perfectly matched discrete layers proposed by
Guddati and Lim [Int. J. Numer. Meth. Engng 66 (6) (2006) 949-977],
including them as type $(L,1)$.
An analysis of the proposed ABCs is performed motivated by the work of
Ainsworth [J. Comput. Phys. 198 (1) (2004) 106-130].
The new ABCs facilitate numerical implementations
of the Helmholtz problem with ABCs
if $Q_N$ finite elements are used in the physical domain
as well as give more insight into this field for the further advancement.
\end{abstract}
\maketitle

\allowdisplaybreaks

% -- Introduction

\section{Introduction}\label{S:introduction}
The Helmholtz equation is of interest
in many fields in physics and engineering.
We are particularly interested in solving the Helmholtz equation
in unbounded domains.
As the standard numerical methods such as finite element methods
require finite computational domains,
numerically solving such problem requires one to map
the original problem in an unbounded domain
to one in a bounded domain.
This is commonly achieved by truncating the unbounded domain
and applying an artificial boundary condition representing
the original unboundedness
on or in the vicinity of the surface of truncation.
Such artificial boundary conditions are called
\emph{absorbing boundary conditions} (ABCs).
The \emph{accuracy}
of an ABC is often measured by the \emph{reflection coefficient},
the ratio of the spurious incoming wave to the outgoing wave.
Many ABCs have been proposed in the last decades, and,
among them, the complete radiation boundary conditions (CRBCs) and
the perfectly matched discrete layers (PMDLs) are
known to be especially accurate and efficient.

ABCs have typically been studied for the wave equation, from which
those for the Helmholtz equation can readily be derived.
The early ABCs proposed and studied in
\cite{EngquistMajda1977,EngquistMajda1979,BaylissTurkel1980,Higdon1986,Higdon1987}
involved high-order temporal and spatial derivatives
on the artificial boundary.
The Higdon ABCs \cite{Higdon1986} were
designed to annihilate $L$ purely propagating plane waves
and included $L$-th order temporal and normal derivatives
of the solution on the boundary.
The Givoli-Neta ABCs \cite{GivoliNeta2003}
introduced auxiliary functions on the boundary and
recast the classic Higdon ABCs
as systems of $L$ equations that only involved first-order
temporal and normal derivatives,
greatly facilitating the numerical implementations.
Hagstrom-Warburton ABCs \cite{HagstromWarburton2004}
\emph{symmetrised} the Givoli-Neta ABCs,
squaring the reflection coefficients.
The CRBCs \cite{HagstromWarburton2009}
further improved the behaviour by directly taking into account
the waves that decays while propagating.
The significance of treating such wave mode was also noted in \cite{Higdon1987}.
This class of ABCs are generally known to show excellent performances
regardless of the choice of the parameters.
Since it was demanded that the auxiliary functions
were only defined \emph{on} the boundary,
the above classes of ABCs required reformulations to remove
the unfavorable first-order normal derivatives.
In \cite{GivoliNeta2003} the Givoli-Neta ABCs were rewritten in forms
that only involved up to second-order temporal and \emph{tangential} derivatives;
such formulations are thus called \emph{second-order formulations}.
The second-order formulations of the later ABCs were derived in
\cite{HagstromWarburton2004,Modave2017}.
Further, in multi dimensions,
the above classes of ABCs required special treatment on edges and corners
on which multiple faces met.
The \emph{edge/corner compatibility conditions} were derived for
the Hagstrom-Warburton ABCs in
\cite{HagstromWarburton2004}
and for the later ABCs in
\cite{Modave2017}.
The second-order formulations as well as
the edge/corner compatibility conditions
require one to work with complex systems of equations on the artificial boundary
that are distinct from the wave equation or the Helmholtz equation
solved in the physical domain.
As noted shortly, the PMDLs have a clear advantage in this regard
when the Helmholtz equation is considered.
Having been developed for the wave equation, however,
the CRBCs and the other related ABCs typically handle
the wave equation more directly and efficiently, but,
as we focus on the Helmholtz equation in this work,
we leave further discussion on the wave equation for future work.

The perfectly matched layers (PMLs) formulation is another class of
conditions widely used to approximate the effect of the unboundedness.
The PMLs were originally proposed in \cite{Berenger1994}
for the time-domain electromagnetics, in which
artificial damping terms were added in the vicinity of
the domain of interest.
We call the domain of interest the \emph{physical domain},
the domain in which damping terms were added
the \emph{artificial domain}, and
call the union of the two
the \emph{computational domain}.
The addition of the artificial damping was later reinterpreted
as a \emph{complex coordinate transformation} in the frequency-domain
in \cite{ChewWeedon1994}.
As for the Helmholtz equation,
the PMLs have been more popular than the CRBCs and the other related ABCs,
presumably because, a single \emph{Helmholtz-like equation},
which reduces to the standard Helmholtz equation
for a particular choice of parameters,
monolithically governs the entire computational domain
even in existence of edges and corners,
which facilitates numerical implementations;
see, e.g., \cite{Berenger1994} for the initial treatment of the corners.
Despite the popularity of the PMLs, however,
the accuracy of the PMLs are known to be sensitive to
the choice of the parameters, which are hard to tune \emph{a priori},
and the ABCs often outperform the PMLs in terms of accuracy.
An attempt to optimise the performance of the PMLs
for a particular finite difference scheme \cite{Asvadurov2003}
then inspired the development of the PMDLs \cite{Guddati2006,GuddatiLim2006}
designed for finite element methods.
The PMDLs formulation is effectively obtained by
first writing the PMLs formulation in the weak form,
then using the \emph{linear} finite element approximation
with \emph{midpoint} reduced integration in the artificial domain
in the direction normal to the artificial boundary.
It has been shown that the PMDLs are essentially equivalent to the CRBCs, but
the PMDLs being derived from the PMLs,
they inherit the favourable property of the PMLs, and
allow for working with the Helmholtz-like equation
across the computational domain even if edges and corners exist.
Specifically, if one chooses to use the $Q_1$ tensor-product Lagrange finite element
in the physical domain,
that Helmholtz-like equation can be discretised monolithically
on the entire computational domain,
leaving one only having to use different quadrature rules
in the physical and artificial domains.
If one chooses to use the $Q_N$ element with $N>1$ in the physical domain,
the artificial domain must be treated separately
as the PMDLs require linear approximation in the artificial domain
in the direction normal to the boundary.
Even in such cases, however,
the PMDLs are generally more tractable than the CRBCs.

In this work we propose a new class of ABCs
for the Helmholtz equation that generalises the PMDLs.
The proposed ABCs are derived from a certain simple class of
perfectly matched layers with $L$ discrete layers in the artificial domain,
using the $Q_N$ tensor-product Lagrange finite elements in conjunction with
the $N$-point Gauss-Legendre quadrature reduced integration
in the artificial domain.
Thus, the proposed ABCs can be classified by a tuple $(L,N)$,
with the PMDLs represented by $(L,1)$.
The ABCs of type $(L,N)$ and those of type $(LN,1)$
introduce the same number of auxiliary degrees of freedom
in the artificial domain, and they are as accurate as each other.
Like the PMDLs, the proposed ABCs only involve a single
Helmholtz-like equation over the entire computational domain
even if edges and corners exist, but
they further allow for monolithic discretisations
over the computational domain
even when the $Q_N$ finite elements with $N>1$ are used in the physical domain.
We present an analysis of the proposed ABCs,
strongly motivated by the techniques used in \cite{Ainsworth2004},
in which dispersive and dissipative behaviour of
high-order discontinuous Galerkin methods were studied,
and in the related works \cite{Ainsworth2004-1,AinsworthWajid2010,Ainsworth2014},
that is to give more insight into the problems of ABCs.

In Sec.~\ref{S:helmholtz} and Sec.~\ref{S:pml}
we introduce the Helmholtz equation and the standard PMLs formulation.
In Sec.~\ref{S:const} we study a certain simple PMLs formulation in detail,
which will aid us with developing the new ABCs.
In Sec.~\ref{S:weak} we introduce the discretised weak forms of the PMLs formulations
introduced in Sec.~\ref{S:pml} and Sec.~\ref{S:const}.
In Sec.~\ref{S:abc} we introduce and study the new class of ABCs.
In Sec.~\ref{S:example} we show one- and three-dimensional examples.
In Sec.~\ref{S:conclusion} we summarise this work and discuss future works.

\section{The Helmholtz equation}\label{S:helmholtz}
The Helmholtz equation in $\real^d$, $d\in\mathbb{N}$,
to be solved for $u(x_1,x_2,\cdots,x_d)$ is given by:
\begin{align}
s^2u-u_{,ii}&=f\quad\textrm{in }\real^d,
\label{E:helmholtz_xys_f}
\end{align}
where sum is taken over $i\in\{1,\cdots,d\}$,
$s\in\complex$ is a complex number with $\Real{s}\geq 0$, and
$f$ is a spatial function.
\eqref{E:helmholtz_xys_f} can be obtained from the wave equation by taking the Laplace transform in time, in which case $\Real{s}>0$, or by making the time-harmonic assumption
using the $\exp{(\iu\omega t)}$ convention with $\omega\in\real$, in which case $s=\iu\omega$, and
$\Real{s}=0$.
Assuming that $f$ has a compact support
in $x<-\delta$ for some $\delta>0$, where $x:=x_1$ is the first coordinate variable,
we have:
\begin{align}
s^2u-u_{,ii}&=0\quad\textrm{in }x>-\delta.
\label{E:helmholtz_xys}
\end{align}
In this work we mostly work with the following mapped problem of \eqref{E:helmholtz_xys}
to be solved for $u(x)$
(with an abuse of notation):
\begin{align}
\gamma^2 u- u_{,xx}=0\quad\textrm{in }x>-\delta,
\label{E:helmholtz_xgamma}
\end{align}
where $\gamma:=\sqrt{s^2+ k^2}$,
where $k:=\sqrt{k_ik_i}$ (sum on $i\in\{2,\cdots,d\}$),
where $k_i\in\real$.
\eqref{E:helmholtz_xgamma} is obtained by taking the Fourier transform of \eqref{E:helmholtz_xys}
in all the coordinate variables except in $x:=x_1$ following the Laplace transform in time \cite{HagstromWarburton2009},
in which case, $k_i$ is the dual variable of $x_i$, or
by making the plane wave assumption using the $\exp{(\iu k_i x_i)}$ (sum on $i\in\{2,\cdots,d\}$) convention
following the time-harmonic assumption.
In the former case, $\Real{s}>0$, and we consider a Riemann surface with branch cut on
$(-\iu\infty,-\iu k)\cup(+\iu k,+\iu\infty)$ and
choose the branch of $\sqrt{s^2+ k^2}$ so that $\gamma=+k$ when $s=0$.
In this case we have $\Real{\gamma}>0$.
In the latter case, $s=\iu\omega$, and we define the value of $\gamma$ as the limit in the above branch
along $s=\iu\omega+\epsilon$ as $\epsilon\rightarrow 0^+$.
In this case we assume that $s\neq\pm\iu k$.
The following two properties then hold for $\Real{s}\geq 0$:
\begin{subequations}
\begin{align}
&[\Imag{\gamma}\!\cdot\!\Imag{s}>0\textrm{ and }\Real{\gamma}\geq 0]
\text{ or }
[\Imag{\gamma}=0\textrm{ and }\Real{\gamma}>0],
\label{E:gamma_prop_1}\\
&\textrm{Re}\hspace{-2pt}\left(\frac{\gamma}{a\!\cdot\!s+b}\right)>0,\quad\forall a,b>0.\label{E:gamma_prop_2}
\end{align}
\label{E:gamma_prop}%
\end{subequations}

One can readily solve \eqref{E:helmholtz_xgamma}
up to two undetermined coefficients, $c_-$ and $c_+$, as:
\begin{align}
u(x)=c_-e^{+\gamma x}+c_+e^{-\gamma x}.
\label{E:helmholtz_xgamma_solution}
\end{align}
Noting \eqref{E:gamma_prop_1}, one can see that
$e^{+\gamma x}$ represents the wave
that is \emph{incoming} and/or grows as $x\rightarrow+\infty$ and
$e^{-\gamma x}$ represents the wave
that is \emph{outgoing} and/or decays as $x\rightarrow+\infty$.
The pure outgoing/decaying wave would thus require $c_{-}=0$, or:
\begin{align}
c_-=0\!\cdot\!c_+.
\label{E:helmholtz_reflection}
\end{align}
This ratio of $c_-$ to $c_+$ is called the
\emph{reflection coefficient}, and is an important measure
of the accuracy of absorbing boundary conditions (ABCs).
On the other hand, the pure outgoing/decaying wave solution would also satisfy:
\begin{align}
u_{,x}(0)+\gamma u(0)=0,
\label{E:helmholtz_sommerfeld}
\end{align}
which is known as the \emph{Sommerfeld radiation condition} (at $x=0$).
The extent to which a given ABC
approximates \eqref{E:helmholtz_sommerfeld}
gives another important measure of the accuracy of that condition.

Our domain of interest, or the \emph{physical domain}, being $x<0$,
we are to truncate the right half-domain.
Although \eqref{E:helmholtz_sommerfeld} gives the \emph{exact}
boundary condition for the outgoing/decaying wave at $x=0$
for the one-dimensional problem \eqref{E:helmholtz_xgamma}
for a given $k\in\real$,
its multidimensional counterpart that is exact for all $k$ does not exist.
In the following sections
we work with several conditions applied in the vicinity of $x=0$
outside the physical domain, which we call the \emph{artificial domain},
that are to approximate the effect of the unboundedness of $x>0$,
ideally only admitting the outgoing/decaying wave at $x=0$.

\section{Perfectly matched layers}\label{S:pml}
The perfectly matched layers (PMLs) method was originally developed
for the time-domain electromagnetics in \cite{Berenger1994}
by adding artificial damping terms in the artificial domain
so that the outgoing wave attenuates.
The addition of artificial damping was later reinterpreted
as \emph{complex coordinate transformation} in \cite{ChewWeedon1994}.
In this section we briefly introduce the PMLs formulation using the idea of
complex coordinate transformation.
Although the PMLs formulation is conventionally derived for a simplified version
of \eqref{E:helmholtz_xgamma}, which
one can obtain by setting $s=\iu\omega$ and $k=0$, or $\gamma=\iu\omega$,
we here use \eqref{E:helmholtz_xgamma} directly
to take into account the effect of multi dimensions.

The crucial step in deriving the PMLs formulation is to introduce
a \emph{complex coordinate transformation},
$\tilde{\chi}:\real\rightarrow\complex$,
which is sufficiently smooth and
satisfies $\tilde{\chi}(x)=x$ in $x<0$.
We denote the transformed coordinate by $\tilde{x}:=\tilde{\chi}(x)$,
and consider solving the following modified version of \eqref{E:helmholtz_xgamma}
for $\tilde{u}(\tilde{x})$:
\begin{align}
\gamma^2\tilde{u}-\tilde{u}_{,\tilde{x}\tilde{x}}=0\quad x>-\delta.
\label{E:pml_inf_tilde}
\end{align}
Defining $u^{\dagger}:=\tilde{u}\circ\tilde{\chi}$,
one can rewrite \eqref{E:pml_inf_tilde} in terms of $x$
as a problem to be solved for $u^{\dagger}(x)$ as:
\begin{align}
\frac{\gamma^2}{g^{\PML}}u^{\dagger}
-(g^{\PML}u^{\dagger}_{,x})_{,x}&=0\quad x>-\delta,
\label{E:pml_inf}
\end{align}
where:
\begin{align}
g^{\PML}(x):=\frac{\dif x}{\dif\tilde{x}},
\label{E:pml_g}
\end{align}
where the superscript $(\cdot)^{\PML}$ is used to represent PMLs.
Note that $g^{\PML}(x)=1$ in $x<0$, so \eqref{E:pml_inf}
reduces to the standard Helmholtz equation \eqref{E:helmholtz_xgamma} in $x<0$.
By construction, the solution to \eqref{E:pml_inf} is readily obtained as:
\begin{align}
u^{\dagger}(x)=
c^{\dagger}_-e^{+\gamma\tilde{\chi}(x)}+
c^{\dagger}_+e^{-\gamma\tilde{\chi}(x)}.
\label{E:pml_inf_solution}
\end{align}
The complex transformation function $\tilde{\chi}(x)$ is chosen so that
the outgoing/decaying wave, represented by $\exp{(-\gamma\tilde{\chi}(x))}$,
attenuates as $x\rightarrow +\infty$.
A popular choice in the common setting of $\gamma=\iu\omega$ is
$\tilde{\chi}(x)=x+(\omega_0/\iu\omega)\!\cdot\!\int_{-\infty}^xf^{\PML}(\tau)\dif\tau$, where
$\omega_0$ is a positive constant for non-dimensionalisation and
$f^{\PML}(x)=0$ in $x\leq 0$ and
$f^{\PML}(x)>0$ in $x>0$.
The outgoing wave solution then becomes
$\exp{(-\gamma x)}\!\cdot\!\exp{(-\omega_0\int_{-\infty}^xf^{\PML}(\tau)\dif\tau)}$,
which coincides with the exact outgoing wave solution $\exp{(-\gamma x)}$
in $x<0$ and decays exponentially as $x\rightarrow +\infty$.
For convenience, we here introduce a set of grid points,
$\{x_{-1},x_{0},\cdots,x_L\}$, where
$-\delta<x_{-1}<x_0=0<x_1<\cdots<x_L$.
The above observation justifies one in truncating the right half-space $x>x_L$
for $x_L$ sufficiently large such that $u^{\dagger}(x_L)\approx0$.
We thus obtain the PMLs formulation to be solved for $u^{\PML}(x)$
up to one undetermined coefficient as:
\begin{subequations}
\begin{align}
\frac{\gamma^2}{g^{\PML}}u^{\PML}
-(g^{\PML}u^{\PML}_{,x})_{,x}&=0
\quad\textrm{ in }x_{-1}<x<x_L,
\label{E:pml_domain}\\
u^{\PML}(x_L)&=0,
\label{E:pml_end}
\end{align}
\label{E:pml}%
\end{subequations}
\eqref{E:pml_end} is called the \emph{termination condition}.
\eqref{E:pml_domain} can readily be solved up to two undetermined
coefficients, $c^{\PML}_-$ and $c^{\PML}_+$, as:
\begin{align}
u^{\PML}(x)=c^{\PML}_-e^{+\gamma\tilde{\chi}(x)}+c^{\PML}_+e^{-\gamma\tilde{\chi}(x)}.
\label{E:pml_solution}
\end{align}
Applying \eqref{E:pml_end}, one can solve \eqref{E:pml}
up to one undetermined coefficient as desired.
The reflection coefficient of the solution thus obtained is given as:
\begin{align}
c^{\PML}_-=-\left(e^{-\gamma\tilde{\chi}(x_L)}\right)^2\!\cdot\!c^{\PML}_+,
\label{E:pml_reflection}
\end{align}
and the solution satisfies the following approximate Sommerfeld radiation condition:
\begin{align}
u^{\PML}_{,x}(0)+\gamma\left(\frac{1+\left(e^{-\gamma\tilde{\chi}(x_L)}\right)^2}{1-\left(e^{-\gamma\tilde{\chi}(x_L)}\right)^2}\right)u^{\PML}(0)=0.
\label{E:pml_sommerfeld}
\end{align}
In the spatially continuous setting
$\exp{(-\gamma\tilde{\chi}(x_L))}^2$
can be made as small as desired by manipulating $\tilde{\chi}(x)$, and
the expression for the reflection coefficient \eqref{E:pml_reflection} and
the approximate Sommerfeld radiation condition \eqref{E:pml_sommerfeld}
can be made arbitrarily close to the exact expressions
\eqref{E:helmholtz_reflection} and \eqref{E:helmholtz_sommerfeld}.
As mentioned in Sec.~\ref{S:weak}, however, this will no longer be the case
after discretisation.
Finally, we note that,
if $g^{\PML}(x_L)u^{\PML}_{,x}(x_L)+\gamma u^{\PML}(x_L)=0$,
which we call the \emph{transformed Sommerfeld radiation condition},
was applied at $x=x_L$ in place of the termination condition \eqref{E:pml_end},
we would obtain $c^{\PML}_{-}=0$ and
$u^{\PML}(x)$ would represent the exact outgoing/decaying wave solution in $x<0$
despite the existence of the PMLs.

\section{Layerwise constant perfectly matched layers}\label{S:const}
Before introducing the new ABCs,
it is instructive to study the \emph{layerwise constant perfectly matched layers}.
We first define a layerwise constant complex function, $g^{\CPML}(x)$, as:
\begin{align}
g^{\CPML}(x):=\gamma_l\quad\textrm{ in }x_{l-1}<x\leq x_l,
\quad l\in\{0,\cdots,L\},
\label{E:const_g}
\end{align}
where $\gamma_l\in\complex$ is a complex constant and 
the superscript $(\cdot)^{\CPML}$ is used
to represent layerwise constant PMLs.
$\gamma_0$ is defined to be $1$, and
$\gamma_l$, $l\in\{1,\cdots,L\}$, is chosen so that $\Real{\gamma/\gamma_l}>0$;
this is achieved by letting $\gamma_l=a\cdot s + b$ for some $a,b>0$,
according to Property \eqref{E:gamma_prop_2}.

We solve a problem that is similar to \eqref{E:pml}, but
explicitly handling the discontinuities of $g^{\CPML}(x)$.
Specifically, we define $u^{\CPML}(x)$ as:
\begin{align}
u^{\CPML}(x):=u^{\CPML}_l(x)\quad\textrm{ on }[x_{l-1},x_l],
\quad l\in\{0,\cdots,L\},
\end{align}
solving the following system of equations for $u^{\CPML}_l(x)$, $l\in\{0,\cdots,L\}$,
which is assumed to be sufficiently smooth
in $(x_{l-1}-\epsilon,x_l+\epsilon)$ for some $\epsilon>0$:
\begin{subequations}
\begin{align}
\frac{\gamma^2}{\gamma_l}u^{\CPML}_l-(\gamma_lu^{\CPML}_{l,x})_{,x}&=0
\quad\textrm{ in }(x_{l-1},x_l),
&&\hspace{-40pt}\forall l\in\{0,\cdots,L\},
\label{E:const_domain}\\
\gamma_lu^{\CPML}_{l,x}(x_l)-
\gamma_{l+1}u^{\CPML}_{l+1,x}(x_l)&=0,
&&\hspace{-40pt}\forall l\in\{0,\cdots,L-1\},
\label{E:const_interface_ux}\\
u^{\CPML}_l(x_l)-u^{\CPML}_{l+1}(x_l)&=0,
&&\hspace{-40pt}\forall l\in\{0,\cdots,L-1\},
\label{E:const_interface_u}\\
u^{\CPML}_L(x_L)&=0.&&
\label{E:const_end}
\end{align}
\label{E:const}%
\end{subequations}
Once we solve \eqref{E:const_domain} up to two undetermined coefficients
for each $l$,
we will be left with a system of $L+L+1=2(L+1)-1$ equations
for $2(L+1)$ unknowns.

\subsection{Domain equations}\label{SS:const_domain}
One can readily solve \eqref{E:const_domain} for $u^{\CPML}_l(x)$
for each $l\in\{0,\cdots,L\}$ up to two undetermined coefficients,
$c^{\CPML}_{l,\textrm{even}}$ and $c^{\CPML}_{l,\textrm{odd}}$, as:
\begin{align}
u^{\CPML}_l(x)=
c^{\CPML}_{l,\textrm{even}}\!\cdot\!\left(
\frac{e^{+\frac{\gamma}{\gamma_l}\bigl(x-x_{l-\frac{1}{2}}\bigr)}+e^{-\frac{\gamma}{\gamma_l}\bigl(x-x_{l-\frac{1}{2}}\bigr)}}{2}
\right)+
c^{\CPML}_{l,\textrm{odd}}\!\cdot\!\left(
\frac{e^{+\frac{\gamma}{\gamma_l}\bigl(x-x_{l-\frac{1}{2}}\bigr)}-e^{-\frac{\gamma}{\gamma_l}\bigl(x-x_{l-\frac{1}{2}}\bigr)}}{2}
\right),
\label{E:const_domain_solution}
\end{align}
where $x_{l-\frac{1}{2}}=(x_{l-1}+x_l)/2$.
We are now left with $2(L+1)-1$ equations for $2(L+1)$ unknowns.

\subsection{Interface conditions}\label{SS:const_interface}
Using \eqref{E:const_domain_solution},
\eqref{E:const_interface_u} and \eqref{E:const_interface_ux}
can be written as:
\begin{align}
\begin{pmatrix}
S^{\CPML}_{\textrm{even}}(\alpha_l) & S^{\CPML}_{\textrm{odd}}(\alpha_l)\\
S^{\CPML}_{\textrm{odd}}(\alpha_l) & S^{\CPML}_{\textrm{even}}(\alpha_l)
\end{pmatrix}
\hspace{-5pt}
\begin{pmatrix}
c^{\CPML}_{l,\textrm{even}}\\
c^{\CPML}_{l,\textrm{odd}}
\end{pmatrix}
-
\begin{pmatrix}
\phantom{-\hspace{-2pt}}S^{\CPML}_\textrm{even}(\alpha_{l+1}) &
\hspace{-10pt}-\hspace{-2pt}S^{\CPML}_\textrm{odd}(\alpha_{l+1})\\
-\hspace{-2pt}S^{\CPML}_\textrm{odd}(\alpha_{l+1}) &
\hspace{-10pt}\phantom{-\hspace{-2pt}}S^{\CPML}_\textrm{even}(\alpha_{l+1})
\end{pmatrix}
\hspace{-5pt}
\begin{pmatrix}
c^{\CPML}_{l+1,\textrm{even}}\\
c^{\CPML}_{l+1,\textrm{odd}}
\end{pmatrix}
=
\begin{pmatrix}
0\\
0
\end{pmatrix}
,
\label{E:const_interface_1}
\end{align}
where $S^{\CPML}_\textrm{even}(\alpha_l)$ and $S^{\CPML}_\textrm{odd}(\alpha_l)$ are
defined as:
\begin{subequations}
\begin{align}
S^{\CPML}_\textrm{even}(\alpha_l):&=
\frac{e^{+\frac{\alpha_l}{2}}+e^{-\frac{\alpha_l}{2}}}{2},\\
S^{\CPML}_\textrm{odd}(\alpha_l):&=
\frac{e^{+\frac{\alpha_l}{2}}-e^{-\frac{\alpha_l}{2}}}{2},
\end{align}
\label{E:const_S}%
\end{subequations}
where:
\begin{align}
\alpha_l:=\frac{\gamma}{\gamma_l}h_l,
\label{E:alpha_l}
\end{align}
where $h_l:=x_l-x_{l-1}$.
Note that $\gamma_0=1$ and we chose
$\gamma_l$, $l\in\{1,\cdots,L\}$, so that $\Real{\alpha_l}>0$.
Performing eigendecompositions of the matrices appearing
in \eqref{E:const_interface_1}, we have:
\begin{subequations}
\begin{align}
\begin{pmatrix}
S^{\CPML}_{\textrm{even}}(\alpha_l) & S^{\CPML}_{\textrm{odd}}(\alpha_l)\\
S^{\CPML}_{\textrm{odd}}(\alpha_l) & S^{\CPML}_{\textrm{even}}(\alpha_l)
\end{pmatrix}
&=
\Phi
\begin{pmatrix}
e^{+\frac{\alpha_l}{2}} & \hspace{-5pt}0\\
0 & \hspace{-5pt}e^{-\frac{\alpha_l}{2}}
\end{pmatrix}
\Phi^\textrm{T},\\
\begin{pmatrix}
\phantom{-\hspace{-2pt}}S^{\CPML}_{\textrm{even}}(\alpha_{l+1}) &
\hspace{-10pt}-\hspace{-2pt}S^{\CPML}_{\textrm{odd}}(\alpha_{l+1})\\
-\hspace{-2pt}S^{\CPML}_{\textrm{odd}}(\alpha_{l+1}) &
\hspace{-10pt}\phantom{-\hspace{-2pt}}S^{\CPML}_{\textrm{even}}(\alpha_{l+1})
\end{pmatrix}
&=
\Phi
\begin{pmatrix}
e^{-\frac{\alpha_{l+1}}{2}} & \hspace{-5pt}0\\
0 & \hspace{-5pt}e^{+\frac{\alpha_{l+1}}{2}}
\end{pmatrix}
\Phi^\textrm{T},
\end{align}
\label{E:const_eigen}%
\end{subequations}
where $\Phi$ is a $2\times 2$ orthonormal matrix of eigenvectors
that is independent of the layer parameter $l$:
\begin{align}
\Phi=
\begin{pmatrix}
\frac{1}{\sqrt{2}} & \hspace{-10pt}\phantom{-\hspace{-0pt}}\frac{1}{\sqrt{2}}\\
\frac{1}{\sqrt{2}} & \hspace{-10pt}-\hspace{-0pt}\frac{1}{\sqrt{2}}
\end{pmatrix}
.
\label{E:Phi}
\end{align}
Using \eqref{E:const_eigen} and \eqref{E:Phi},
\eqref{E:const_interface_1} reduces to:
\begin{align}
\begin{pmatrix}
\frac{e^{+\frac{\alpha_l}{2}}}{e^{-\frac{\alpha_{l+1}}{2}}} & \hspace{-5pt}0\\
0 & \hspace{-5pt}\frac{e^{-\frac{\alpha_l}{2}}}{e^{+\frac{\alpha_{l+1}}{2}}}
\end{pmatrix}
\hspace{-5pt}
\begin{pmatrix}
c^{\CPML}_{l,-}\\
c^{\CPML}_{l,+}
\end{pmatrix}
=
\begin{pmatrix}
c^{\CPML}_{l+1,-}\\
c^{\CPML}_{l+1,+}
\end{pmatrix}
,
\label{E:const_recursion}
\end{align}
where:
\begin{align}
\begin{pmatrix}
c^{\CPML}_{l,-}\\
c^{\CPML}_{l,+}
\end{pmatrix}
=
\Phi^\textrm{T}
\begin{pmatrix}
c^{\CPML}_{l,\textrm{even}}\\
c^{\CPML}_{l,\textrm{odd}}
\end{pmatrix}
.
\label{E:const_cl+-}
\end{align}
Using \eqref{E:const_recursion}, one can recursively eliminate
$c^{\CPML}_{l,-}$ and $c^{\CPML}_{l,+}$, $l\in\{1,\cdots,L-1\}$, to obtain:
\begin{align}
\begin{pmatrix}
\prod\limits_{l=0}^{L-1}\frac{e^{+\frac{\alpha_l}{2}}}{e^{-\frac{\alpha_{l+1}}{2}}} & \hspace{-5pt}0\\
0 & \hspace{-5pt}\prod\limits_{l=0}^{L-1}\frac{e^{-\frac{\alpha_l}{2}}}{e^{+\frac{\alpha_{l+1}}{2}}}
\end{pmatrix}
\hspace{-5pt}
\begin{pmatrix}
c^{\CPML}_{0,-}\\
c^{\CPML}_{0,+}
\end{pmatrix}
=
\begin{pmatrix}
c^{\CPML}_{L,-}\\
c^{\CPML}_{L,+}
\end{pmatrix}
.
\label{E:const_c0_cL}
\end{align}
We are now left with $1$ equation for $2$ unknowns.

$c^{\CPML}_{0,-}$ and $c^{\CPML}_{0,+}$ respectively represent
the incoming/growing wave and the outgoing/decaying wave.
Here, we suppose that there was an additional layer $[x_{-2},x_{-1}]$ with
$-\delta<x_{-2}<x_{-1}<0$ and
define $h_{-1}:=x_{-1}-x_{-2}=h_0$ and $\gamma_{-1}:=1=\gamma_0$.
If we extended the above analysis to include $[x_{-2},x_{-1}]$,
the solution $u_{-1}^{\CPML}$ in $(x_{-2}-\epsilon,x_{-1}+\epsilon)$
would be represented in terms of $c^{\CPML}_{-1,-}$ and $c^{\CPML}_{-1,+}$.
Using \eqref{E:const_recursion} with $l=-1$ and \eqref{E:alpha_l},
$c^{\CPML}_{-1,-}$ and $c^{\CPML}_{-1,+}$ would be related to $c^{\CPML}_{0,-}$ and $c^{\CPML}_{0,+}$ as:
\begin{align}
\begin{pmatrix}
e^{+\gamma h_0} & \hspace{-10pt}0\\
0               & \hspace{-10pt}e^{-\gamma h_0}
\end{pmatrix}
\hspace{-5pt}
\begin{pmatrix}
c^{\CPML}_{-1,-}\\
c^{\CPML}_{-1,+}
\end{pmatrix}
=
\begin{pmatrix}
c^{\CPML}_{0,-}\\
c^{\CPML}_{0,+}
\end{pmatrix}
,
\end{align}
which implies the above statement.

\subsection{Termination condition}\label{SS:const_end}
Finally, we apply the termination condition \eqref{E:const_end}.
Using \eqref{E:const_domain_solution} with $l=L$ along with \eqref{E:const_S},
\eqref{E:const_end} can be written as:
\begin{align}
S^{\CPML}_{\textrm{even}}(\alpha_L)\!\cdot\!c^{\CPML}_{L,\textrm{even}}
+ S^{\CPML}_{\textrm{odd}}(\alpha_L)\!\cdot\!c^{\CPML}_{L,\textrm{odd}}
=0,
\label{E:const_end_u_1}
\end{align}
or, rewriting in terms of $c^{\CPML}_{L,-}$ and $c^{\CPML}_{L,+}$
using \eqref{E:const_cl+-}:
\begin{align}
\frac{1}{\sqrt{2}}
\begin{pmatrix}
e^{+\frac{\alpha_L}{2}}\hspace{5pt}e^{-\frac{\alpha_L}{2}}
\end{pmatrix}
\hspace{-3pt}
\begin{pmatrix}
c^{\CPML}_{L,-}\\
c^{\CPML}_{L,+}
\end{pmatrix}
=0.
\label{E:const_end_u_2}
\end{align}
Using \eqref{E:const_c0_cL} in conjunction with \eqref{E:const_end_u_2}, we obtain:
\begin{align}
c^{\CPML}_{0,-}=
-
\left(
e^{-\alpha_0}
\right)
\prod_{l=1}^{L}
\left(
e^{-\alpha_l}
\right)^2
\!\cdot\!c^{\CPML}_{0,+}
,
\label{E:const_reflection}
\end{align}
which compares with the expression for the reflection coefficient for
the standard PMLs formulation \eqref{E:pml_reflection}.
As it is the case for the standard PMLs,
the reflection coefficient for the layerwise constant PMLs can be made arbitrarily small
by manipulating the values of $\gamma_l$, $l\in\{1,\cdots,L\}$.
in the spatially continuous setting.
The system \eqref{E:const} has now been solved up to one undermined coefficient
as desired.

We now study the ability of the layerwise constant PMLs to approximate
the Sommerfeld radiation condition \eqref{E:helmholtz_sommerfeld}.
For this analysis, \eqref{E:const_domain} only needs to be satisfied
for $l\in\{1,\cdots,L\}$.
Following the same line as the the derivation of \eqref{E:const_reflection},
we have:
\begin{align}
c^{\CPML}_{1,-}=
-
\left(
e^{-\alpha_1}
\right)
\prod_{l=2}^{L}
\left(
e^{-\alpha_l}
\right)^2
\!\cdot\!c^{\CPML}_{1,+}
.
\label{E:const_reflection_1}
\end{align}
Using \eqref{E:const_domain_solution} with $l=1$ along with \eqref{E:const_S},
\eqref{E:const_interface_ux} and \eqref{E:const_interface_u} with $l=0$,
and $\gamma_0=1$, we have:
\begin{subequations}
\begin{align}
u^{\CPML}_0(0)&=
S^{\CPML}_{\textrm{even}}(\alpha_1)\!\cdot\!c^{\CPML}_{1,\textrm{even}}-
S^{\CPML}_{\textrm{odd}}(\alpha_1)\!\cdot\!c^{\CPML}_{1,\textrm{odd}},\\
u^{\CPML}_{0,x}(0)&=
-\gamma\left(
 S^{\CPML}_{\textrm{odd}}(\alpha_1)\!\cdot\!c^{\CPML}_{1,\textrm{even}}
-S^{\CPML}_{\textrm{even}}(\alpha_1)\!\cdot\!c^{\CPML}_{1,\textrm{odd}}
\right)
\end{align}
\label{E:const_approx_1}%
\end{subequations}
or, writing in terms of $c^{\CPML}_{1,-}$ and $c^{\CPML}_{1,+}$
using \eqref{E:const_cl+-}:
\begin{subequations}
\begin{align}
u^{\CPML}_0(0)&=
\frac{1}{\sqrt{2}}
\begin{pmatrix}
e^{-\frac{\alpha_1}{2}}\hspace{5pt}e^{+\frac{\alpha_1}{2}}
\end{pmatrix}
\hspace{-3pt}
\begin{pmatrix}
c^{\CPML}_{1,-}\\
c^{\CPML}_{1,+}
\end{pmatrix}
,\\
u^{\CPML}_{0,x}(0)&=
-\frac{1}{\sqrt{2}}\gamma
\begin{pmatrix}
-e^{-\frac{\alpha_1}{2}}\hspace{5pt}e^{+\frac{\alpha_1}{2}}
\end{pmatrix}
\hspace{-3pt}
\begin{pmatrix}
c^{\CPML}_{1,-}\\
c^{\CPML}_{1,+}
\end{pmatrix}
.
\end{align}
\label{E:const_approx_2}%
\end{subequations}
Using \eqref{E:const_approx_2} in conjunction with
\eqref{E:const_reflection_1}, we obtain
the approximate Sommerfeld radiation condition provided by
the layerwise constant PMLs as:
\begin{align}
u^{\CPML}_{0,x}(0)
+
\gamma
\left(
\frac{1+
\prod\limits_{l=1}^{L}
\left(
e^{-\alpha_l}
\right)^2
}{1-
\prod\limits_{l=1}^{L}
\left(
e^{-\alpha_l}
\right)^2
}
\right)
\hspace{-2pt}
u^{\CPML}_0(0)
=0,
\label{E:const_sommerfeld}
\end{align}
which compares with \eqref{E:pml_sommerfeld}.
Note, in the spatially continuous setting,
the approximation \eqref{E:const_sommerfeld} can be made arbitrarily accurate
by manipulating the values of $\gamma_l$, $l\in\{1,\cdots,L\}$.

\subsection{Transformed Sommerfeld radiation condition}\label{SS:const_transparency}
In this section we consider applying 
the transformed Sommerfeld radiation condition,
$\gamma_Lu^{\CPML}_{L,x}(x_L)+\gamma u^{\CPML}_L(x_L)=0$, at $x=x_L$
in stead of the termination condition \eqref{E:const_end}.
Using \eqref{E:const_domain_solution} with $l=L$
along with \eqref{E:const_S},
the transformed Sommerfeld radiation condition can be written as:
\begin{align}
\gamma\left(
S^{\CPML}_{\textrm{odd}}(\alpha_L)\!\cdot\!c^{\CPML}_{L,\textrm{even}}
+ S^{\CPML}_{\textrm{even}}(\alpha_L)\!\cdot\!c^{\CPML}_{L,\textrm{odd}}
\right)
+
\gamma\left(
S^{\CPML}_{\textrm{even}}(\alpha_L)\!\cdot\!c^{\CPML}_{L,\textrm{even}}+
S^{\CPML}_{\textrm{odd}}(\alpha_L)\!\cdot\!c^{\CPML}_{L,\textrm{odd}}
\right)
=0.
\label{E:const_end_ux_1}
\end{align}
Simplifying the expression and rewriting in terms of $c^{\CPML}_{L,-}$ and $c^{\CPML}_{L,+}$
using \eqref{E:const_cl+-}, we have:
\begin{align}
\sqrt{2}\gamma
\begin{pmatrix}
e^{+\frac{\alpha_L}{2}}\hspace{5pt}0
\end{pmatrix}
\hspace{-3pt}
\begin{pmatrix}
c^{\CPML}_{L,-}\\
c^{\CPML}_{L,+}
\end{pmatrix}
=0.
\label{E:const_end_ux_2}
\end{align}
Using \eqref{E:const_c0_cL} in conjunction with \eqref{E:const_end_ux_2},
we obtain:
\begin{align}
c^{\CPML}_{0,-}=0\!\cdot\!c^{\CPML}_{0,+},
\label{E:const_sommerfeld_transformed_reflection}
\end{align}
which implies that the layerwise constant PMLs with
the transformed Sommerfeld radiation condition at $x=x_L$
only admit the outgoing/decaying wave at $x=0$.

Following the same lines as the derivation of \eqref{E:const_sommerfeld},
one can show that the layerwise constant PMLs with
the transformed Sommerfeld radiation condition at $x=x_L$
provide the exact Sommerfeld radiation condition at $x=0$
for the problem solved in the physical domain, i.e.,:
\begin{align}
u^{\CPML}_{0,x}(0)+\gamma u^{\CPML}_0(0)=0.
\label{E:const_sommerfeld_transformed_sommerfeld}
\end{align}

\section{Weak form and finite element analysis}\label{S:weak}
We let $D:=(x_{-1},x_L)$, and define a mesh $\mathcal{T}$ of $D$ as
$\mathcal{T}:=\{[x_{l-1},x_l]\}_{l\in\{0,\cdots,L\}}$.
We let $V^{\FE}$ be a function space on $\mathcal{T}$ defined as:
\begin{align}
V^{\FE}:=\{w\in H^1(D)|w_{|K}\in\mathbb{P}_N,K\in\mathcal{T}\},
\end{align}
where $\mathbb{P}_N$ is the vector space of polynomials of degree at most $N$.
$V^{\FE}$ can be decomposed into a direct sum of subspaces as:
\begin{align}
V^{\FE}=V^{\FE}_\sharp\oplus V^{\FE}_\flat,
\end{align}
where:
\begin{subequations}
\begin{align}
V^{\FE}_\sharp&:=\textrm{span}(\{w_{-1}^{\HAT},w_0^{\HAT},\cdots,w_{L}^{\HAT}\}),\\
V^{\FE}_\flat&:=\{w\in V^{\FE}|w(x_l)=0,\quad\forall l\in\{-1,0,\cdots,L\}\},
\end{align}
\label{E:direct_sum}%
\end{subequations}
where $w_l^{\HAT}$ represents \emph{hat} functions defined as:
\begin{subequations}
\begin{align}
w_{-1}^{\HAT}&:=\left\{
\begin{array}{cl}
(x_0-x)/(x_0-x_{-1}), &x_{-1}\leq x<x_0\\
0, &\textrm{otherwise}
\end{array}
\right.,\\
w_l^{\HAT}&:=\left\{
\begin{array}{cl}
(x-{x_{l-1}})/(x_l-x_{l-1}), &x_{l-1}\leq x<x_l\\
(x_{l+1}-x)/(x_{l+1}-x_l), &x_{l\phantom{-1}}\leq x<x_{l+1},
\quad\forall l\in\{0,\cdots,L-1\}\\
0, &\textrm{otherwise}
\end{array}
\right.,\\
w_L^{\HAT}&:=\left\{
\begin{array}{cl}
(x-{x_{L-1}})/(x_L-x_{L-1}), &x_{L-1}\leq x<x_L\\
0, &\textrm{otherwise}
\end{array}
\right..
\end{align}
\end{subequations}
The weak form of the PMLs formulation \eqref{E:pml}
to be used with the finite element method
that is solved for $u^{\FE}\in V^{\FE}$ is characterised by:
\begin{subequations}
\begin{align}
\int_{x_{-1}}^{x_L}\frac{\gamma^2}{g^{\PML}}u^{\FE}\ol{w^{\FE}}
+g^{\PML}u^{\FE}_{,x}\ol{w^{\FE}_{,x}}\dif x&=0,
\quad\forall w^{\FE}\in V^{\FE}\backslash\{w_{-1}^{\HAT},w_L^{\HAT}\},
\label{E:pml_weak_domain}\\
{u^{\FE}}(x_{L})&=0.
\label{E:pml_weak_end}
\end{align}
\label{E:pml_weak}%
\end{subequations}
\eqref{E:pml_weak} can be solved up to one undetermined coefficient.
The layerwise constant PMLs \eqref{E:const} are also characterised
by the same weak form \eqref{E:pml_weak}
with $g^{\PML}(x)$ \eqref{E:pml_g} replaced by $g^{\CPML}(x)$ \eqref{E:const_g},
but an alternative weak form directly derived from \eqref{E:const}
will turn out to be useful in our analysis presented in Sec.~\ref{S:abc}.
Specifically, we define $u^{\FE}(x)$ as:
\begin{align}
u^{\FE}(x):=u^{\FE}_l(x)\textrm{ on }[x_{l-1},x_l],\quad l\in\{0,\cdots,L\},
\end{align}
solving the following system
for $u^{\FE}_l(x)$, $l\in\{0,\cdots,L\}$,
a polynomial of order $N$ defined on $[x_{l-1},x_l]$,
up to one undetermined coefficient:
\begin{subequations}
\begin{align}
\int_{x_{l-1}}^{x_l}\frac{\gamma^2}{\gamma_l}u_{l}^{\FE}\ol{w^{\FE}}
+\gamma_lu_{l,x}^{\FE}\ol{w_{,x}^{\FE}}\dif x&=0,\quad\forall w^{\FE}\in V^{\FE}_\flat,\notag\\
&\hspace{10pt}\forall l\in\{0,\cdots,L\},
\label{E:const_weak_domain}\\
\int_{x_{l-1}}^{x_l}\frac{\gamma^2}{\gamma_l}u_l^{\FE}\ol{w_l^{\HAT}}
+\gamma_lu_{l,x}^{\FE}\ol{w_{l,x}^{\HAT}}\dif x+
\int_{x_l}^{x_{l+1}}\frac{\gamma^2}{\gamma_{l+1}}u_{l+1}^{\FE}\ol{w_{l}^{\HAT}}
+\gamma_{l+1}u_{l+1,x}^{\FE}\ol{w_{l,x}^{\HAT}}\dif x&=0,\notag\\
&\hspace{10pt}\forall l\in\{0,\cdots,L- 1\},
\label{E:const_weak_interface_gradu}\\
{u_l^{\FE}}(x_l)-u_{l+1}^{\FE}(x_l)&=0,\notag\\
&\hspace{10pt}\forall l\in\{0,\cdots,L-1\},
\label{E:const_weak_interface_u}\\
{u_L^{\FE}}(x_L)&=0.
\label{E:const_weak_end}
\end{align}
\label{E:const_weak}%
\end{subequations}

Note that, as $g^{\PML}(x)=g^{\CPML}(x)=1$ in $x<0$ by definition,
the integrand of the weak form \eqref{E:pml_weak_domain} reduces to that
of the standard Helmholtz equation in $x<0$,
which greatly simplifies the numerical implementation.
In multidimensional problems edges and corners can be treated
virtually with no added cost.
However, though, in the continuous setting,
one can manipulate $g^{\PML}(x)$ or $g^{\CPML}(x)$
to make the error due to the termination condition,
called the \emph{termination error},
as small as desired,
making the reflection coefficient,
defined in \eqref{E:pml_reflection} or \eqref{E:const_reflection},
arbitrarily small and the approximation,
\eqref{E:pml_sommerfeld} or \eqref{E:const_sommerfeld},
arbitrarily close to the exact Sommerfeld radiation condition,
in the discrete setting represented by \eqref{E:pml_weak},
it is well known that the PMLs suffer from another type of error,
called the \emph{discretisation error},
which is sensitive to the choice of the parameters and
hard to predict \emph{a priori}.
Typically, for a fixed computational cost,
the discretisation error increases when one attempts to decrease
the termination error, and
there is a trade-off between these two types of error;
see, e.g., \cite{BindelGovindjee2005,Koyama2008,Sagiyama2013}.
The CRBCs \cite{HagstromWarburton2009},
and other related ABCs \cite{GivoliNeta2003,HagstromWarburton2004},
are generally more accurate, and their performance is
much less sensitive to the choice of the parameters.
Those ABCs, however, require special treatment on the artificial boundary
\cite{GivoliNeta2003,HagstromWarburton2004,Modave2017}, and implementation
becomes even more involved if edges and corners exist
\cite{HagstromWarburton2004,Modave2017}.
The PMDLs \cite{Guddati2006,GuddatiLim2006} are virtually equivalent to
the CRBCs, but they are more implementation friendly.
The PMDLs formulation is effectively obtained from \eqref{E:pml_weak}
with $g^{\PML}(x)$ replaced by $g^{\CPML}(x)$,
using the $Q_1$ finite element instead of $Q_N$ in the artificial domain, $x>0$,
and using the \emph{midpoint} reduced integration
when evaluating the integrand of \eqref{E:pml_weak_domain}
on each cell in the artificial domain.
Thus, the PMDLs inherit the favourable property of the PMLs, i.e.,
the low cost of implementation, and are attractive alternative to the CRBCs,
potentially except that one is always required to use the $Q_1$ element
in the artificial domain even if the $Q_N$ element with $N>1$ is used
in the physical domain.
In Sec.~\ref{S:abc} we introduce a new class of
absorbing boundary conditions that generalises the PMDLs
in a manner that allows for using $Q_N$ finite elements in the artificial domain.

\section{Gauss-Legendre Quadrature Reduced-Integration absorbing boundary conditions}\label{S:abc}
Inspired by the PMDLs formulation \cite{Guddati2006,GuddatiLim2006},
we use the weak form \eqref{E:pml_weak}
with $g^{\PML}(x)$ \eqref{E:pml_g} replaced
with $g^{\CPML}(x)$ \eqref{E:const_g}
that have $L$ absorbing layers in the artificial domain,
i.e., $\{[x_{l-1},x_l]\}_{l\in\{1,\cdots,L\}}$,
and use the $Q_N$ finite element
in the artificial domain, $x>0$, as well as
in the physical domain, $x<0$,
in conjunction with the $N$-point Gauss-Legendre quadrature reduced integration
for evaluation of the integral in \eqref{E:pml_weak_domain}
on each cell in the artificial domain.
The proposed ABCs can thus be classified by a tuple $(L,N)$.
We find that the new ABCs generalise the PMDLs,
including them as type $(L,1)$.
Our analysis is motivated by the technique
introduced in \cite{Ainsworth2004}.
As it is more suitable for the analysis, we employ the equivalent weak form
\eqref{E:const_weak} instead of
directly working with \eqref{E:pml_weak};
in actual numerical implementations,
we use \eqref{E:pml_weak}, or its multidimensional counterparts,
in conjunction with the reduced integration in the artificial domain.
Specifically, we define the solution $u^{\ABC}\in V^{\FE}$ as:
\begin{align}
u^{\ABC}(x):=u_l^{\ABC}(x)\textrm{ on }[x_{l-1},x_l],\quad l\in\{0,\cdots,L\},
\label{E:abc_u}
\end{align}
solving the following system for $u_l^{\ABC}(x)$, $l\in\{0,\cdots,L\}$,
a polynomial of order $N$ defined on $[x_{l-1},x_l]$,
up to one undetermined coefficient:
\begin{subequations}
\begin{align}
\rint_{x_{l-1}}^{x_l}\frac{\gamma^2}{\gamma_l}u_{l}^{\ABC}\ol{w^{\ABC}}
+\gamma_lu_{l,x}^{\ABC}\ol{w_{,x}^{\ABC}}\dif x&=0,
\quad\forall w^{\ABC}\in V^{\FE}_\flat,\notag\\
&\hspace{0pt}\forall l\in\{0,\cdots,L\},
\label{E:abc_domain}\\
\rint_{x_{l-1}}^{x_l}\frac{\gamma^2}{\gamma_l}u_l^{\ABC}\ol{w_l^{\HAT}}
+\gamma_lu_{l,x}^{\ABC}\ol{w_{l,x}^{\HAT}}\dif x+
\rint_{x_l}^{x_{l+1}}\frac{\gamma^2}{\gamma_{l+1}}u_{l+1}^{\ABC}\ol{w_{l}^{\HAT}}
+\gamma_{l+1}u_{l+1,x}^{\ABC}\ol{w_{l,x}^{\HAT}}\dif x&=0,\notag\\
&\hspace{0pt}\forall l\in\{0,\cdots,L- 1\},
\label{E:abc_interface_gradu}\\
u_l^{\ABC}(x_l)-u_{l+1}^{\ABC}(x_l)&=0,\notag\\
&\hspace{0pt}\forall l\in\{0,\cdots,L-1\},
\label{E:abc_interface_u}\\
u_L^{\ABC}(x_L)&=0,
\label{E:abc_end_u}
\end{align}
\label{E:abc}%
\end{subequations}
where $\rint$ represents the $N$-point Gauss-Legendre quadrature
reduced-order integration.
Note that $\gamma_0=1$ and $\gamma_l$, $l\in\{1,\cdots,L\}$,
have been chosen so that $\Real{\alpha_l}>0$,
where $\alpha_l$ is defined in \eqref{E:alpha_l}.
Although \eqref{E:abc_domain} with $l=0$ and
\eqref{E:abc_interface_gradu} with $l=0$
imply that the $N$-point reduced integration is also used in the physical domain,
$x<0$, this is not a requirement for the proposed ABCs to work;
see Sec.~\ref{SS:abc_end}.
\eqref{E:abc} is a system of
$(L+1)(N-1)+L+L+1 = (L+1)(N+1)-1$ equations for
$(L+1)(N+1)$ unknowns.
We let the superscript $(\cdot)^{\ABC}$ represent the ABCs of type $(L,N)$.

For each $l\in\{0,\dots,L\}$, we define the affine geometric mapping
$\bs{T}_l:[-1,1]\rightarrow [x_{l-1},x_l]$,
from the reference interval $[-1,1]$
to the mesh cell $[x_{l-1},x_l]\in\mathcal{T}$,
as $\bs{T}_l(\zeta)=((1-\zeta)x_{l-1}+(1+\zeta)x_l)/2$,
where $\zeta\in [-1,1]$ represents the reference coordinate,
and define
$\uref_l^{\ABC}:=u_l^{\ABC}\circ\bs{T}_l$
as the pullback of $u_l^{\ABC}$ by the geometric mapping $\bs{T}_l$.
We then denote by $P_n^{(p,q)}$ the Jacobi polynomial of type $(p,q)$ and order $n$
(cf. (8.960)\textsubscript{1} of \cite{Gradshtein2014}),
and represent $\uref_l^{\ABC}(\zeta)$
using the Jacobi polynomials as:
\begin{subequations}
\begin{align}
\uref_l^{\ABC}(\zeta)&=
\phantom{\left(\frac{\alpha_l}{2}\right)^2}
\sum_{n=0}^{N}c_{l,n}^{\ABC}\alpha_l^n\frac{(2N- n)!}{(2N)!}P_n^{(N- n,N- n)}(\zeta),
\label{E:ul}\\
\intertext{which, using (8.961)\textsubscript{4} of \cite{Gradshtein2014},
lets us write
$\uref^{\ABC}_{l,\zeta}(\zeta)$ and
$\uref^{\ABC}_{l,\zeta\zeta}(\zeta)$ as:}
\uref^{\ABC}_{l,\zeta}(\zeta)&=
\left(\frac{\alpha_l}{2}\right)^{\phantom{2}}
\sum_{n=0}^{N-1}c_{l,n+1}^{\ABC}\alpha_l^n\frac{(2N- n)!}{(2N)!}P_n^{(N- n,N- n)}(\zeta),
\label{E:ulz}\\
\uref^{\ABC}_{l,\zeta\zeta}(\zeta)&=
\left(\frac{\alpha_l}{2}\right)^2
\sum_{n=0}^{N-2}c_{l,n+2}^{\ABC}\alpha_l^n\frac{(2N- n)!}{(2N)!}P_n^{(N- n,N- n)}(\zeta),
\label{E:ulzz}
\end{align}
\end{subequations}
where $\alpha_l$ is defined in \eqref{E:alpha_l} and $c^{\ABC}_{l,n}$ represents
the constants to be determined.

\subsection{Domain equations}\label{SS:abc_domain}
For each $l\in\{0,\cdots,L\}$,
we solve \eqref{E:abc_domain} up to two undetermined coefficients.
We rewrite \eqref{E:abc_domain} in terms of $\zeta$ as:
\begin{align}
\left(\frac{2}{\alpha_l}\right)\gamma
\rint_{-1}^{1}\left(\frac{\alpha_l}{2}\right)^2\uref_l^{\ABC}\ol{\wref^{\ABC}_l}
+\uref_{l,\zeta}^{\ABC}\ol{\wref_{l,\zeta}^{\ABC}}\dif\zeta=0,
\quad\wref^{\ABC}_l={w^{\ABC}}_{|[x_{l-1,x_l}]}\circ\bs{T}_l,
\quad\forall w^{\ABC}\in V^{\FE}_\flat.
\label{E:abc_domain_1}
\end{align}

We first consider the case in which $N\geq 2$.
As $w^{\ABC}\in V^{\FE}_\flat$, $\wref^{\ABC}_l(\zeta)$ can be written as
$\wref^{\ABC}_l(\zeta)=(1-\zeta)(1+\zeta)\check{w}^{\ABC}(\zeta)$,
where $\check{w}^{\ABC}\in\mathbb{P}_{N-2}$.
Noting that the gradient term is exactly evaluated
even under the $N$-point reduced integration,
one can apply the divergence theorem and rewrite \eqref{E:abc_domain_1} as:
\begin{align}
\left(\frac{2}{\alpha_l}\right)\gamma
\rint_{-1}^{1}\left(
\left(\frac{\alpha_l}{2}\right)^2\uref_l^{\ABC}
-\uref_{l,\zeta\zeta}^{\ABC}
\right)\!\cdot\!(1-\zeta)(1+\zeta)\ol{\check{w}^{\ABC}}\dif\zeta=0,
\quad\forall\check{w}^{\ABC}\in\mathbb{P}_{N-2}.
\label{E:abc_domain_2}
\end{align}
Using \eqref{E:ul} and \eqref{E:ulzz}, \eqref{E:abc_domain_2} can be rewritten as:
\begin{align}
\left(\frac{\alpha_l}{2}\right)\gamma\int_{-1}^{1}
\sum_{n=0}^{N-2}(c_{l,n}^{\ABC}-c_{l,n+2}^{\ABC})\alpha_l^n\frac{(2N- n)!}{(2N)!}P_n^{(N- n,N- n)}
\!\cdot\!(1-\zeta)(1+\zeta)\ol{\check{w}^{\ABC}}
\dif \zeta&\notag\\
+\left(\frac{\alpha_l}{2}\right)\gamma\int_{-1}^{1}
c_{l,N-1}^{\ABC}\alpha_l^{N-1}\frac{(N+ 1)!}{(2N)!}P_{N- 1}^{(1,1)}
\!\cdot\!(1-\zeta)(1+\zeta)\ol{\check{w}^{\ABC}}
\dif\zeta&\notag\\
+\left(\frac{\alpha_l}{2}\right)\gamma\rint_{-1}^{1}
c_{l,N}^{\ABC}\alpha_l^{N}\frac{(N)!}{(2N)!}P_N^{(0,0)}
\!\cdot\!(1-\zeta)(1+\zeta)\ol{\check{w}^{\ABC}}
\dif \zeta&=0,
\quad\forall\check{w}^{\ABC}\in\mathbb{P}_{N-2}.
\label{E:abc_domain_3}
\end{align}
The first and the second terms are evaluated exactly with the $N$-point rule.
The second term vanishes
due to the orthogonality property of the Jacobi polynomial $P_{N-1}^{(1,1)}$;
see (7.391)\textsubscript{4} of \cite{Gradshtein2014}.
The last term vanishes as $N$ Gauss-Legendre quadrature points
are the zeros of $P_N^{(0,0)}$.
We are now left with the first term.
If we write $P_n^{(N- n,N- n)}(\zeta)$, $n\in\{0,\cdots,N-2\}$,
and $\check{w}(\zeta)$ as linear combinations of
$P_{n'}^{(1,1)}(\zeta)$, $n'\in\{0,\cdots,N- 2\}$,
we can write \eqref{E:abc_domain_3} in a matrix form as:
\begin{align}
\left(\frac{\alpha_l}{2}\right)\gamma\int_{-1}^{1}
(1+ \zeta)(1- \zeta)\mathbf{Q}(\zeta)\mathbf{P}(\zeta)^\textrm{T}\mathbf{C}\dif \zeta=
\left(\frac{\alpha_l}{2}\right)\gamma
\left[\int_{-1}^{1}
(1+ \zeta)(1- \zeta)\mathbf{Q}(\zeta)\mathbf{Q}(\zeta)^\textrm{T}\dif \zeta
\right]
\mathbf{A}^\textrm{T}\mathbf{C}
=0,
\label{E:abc_domain_4}
\end{align}
where $\mathbf{Q}(\zeta)$, $\mathbf{P}(\zeta)$, and $\mathbf{C}$
are vectors of shape $(N- 1)\times 1$ defined as:
\begin{align}
\mathbf{Q}(\zeta)&:=[P_{n'}^{(1,1)}(\zeta)]^\textrm{T}_{0\leq n'\leq N-2},\\
\mathbf{P}(\zeta)&:=[P_{n}^{(N-n,N-n)}(\zeta)]^\textrm{T}_{0\leq n\leq N-2}=\mathbf{A}\mathbf{Q}(\zeta),\\
\mathbf{C}\phantom{(\zeta)}&:=[C_{n}]^\textrm{T}_{0\leq n\leq N-2},
\end{align}
where:
\begin{align}
C_n:=(c_{l,n}^{\ABC}-c_{l,n+2}^{\ABC})\alpha_l^n\frac{(2N- n)!}{(2N)!},
\end{align}
and $\mathbf{A}$ is the invertible constant matrix of shape $(N-1)\times(N-1)$
representing the basis transformation.
Due to the orthogonality of $P_n^{(1,1)}$
((7.391)\textsubscript{4} of \cite{Gradshtein2014})
and the invertibility of $\mathbf{A}$,
we have $\mathbf{C}=\mathbf{0}$, or:
\begin{align}
c_{l,n+ 2}^{\ABC}=c_{l,n}^{\ABC},\quad n\in\{0,\cdots,N- 2\}.
\label{E:abc_c=c}
\end{align}
Using \eqref{E:ul} in conjunction with \eqref{E:abc_c=c},
$\uref^{\ABC}_l(\zeta)$ can now be represented by
two undetermined coefficients,
$c_{l,\textrm{even}}^{\ABC}(=c_{l,0}^{\ABC})$ and
$c_{l,\textrm{odd}}^{\ABC}(=c_{l,1}^{\ABC})$
for each $l\in\{0,\cdots,L\}$.

For the case of $N=1$, we simply set
$c_{l,\textrm{even}}^{\ABC}=c_{l,0}^{\ABC}$ and $c_{l,\textrm{odd}}^{\ABC}=c_{l,1}^{\ABC}$.

\subsection{Interface conditions}\label{SS:abc_interface}
For each $l\in\{0,\cdots,L-1\}$,
we write \eqref{E:abc_interface_gradu} in terms of $\zeta$ and obtain:
\begin{align}
\left(\frac{2}{\alpha_l}\right)\gamma
\rint_{-1}^{1}\left(\frac{\alpha_l}{2}\right)^2
\uref_l^{\ABC}\!\cdot\!\left(\frac{1+ \zeta}{2}\right)
+\uref_{l,\zeta}^{\ABC}\!\cdot\!\left(\frac{1+ \zeta}{2}\right)\!\!{}_{,\zeta}
\dif\zeta&\notag\\
+\left(\frac{2}{\alpha_{l+1}}\right)\gamma
\rint_{-1}^{1}\left(\frac{\alpha_{l+1}}{2}\right)^2
\uref_{l+1}^{\ABC}\!\cdot\!\left(\frac{1- \zeta}{2}\right)
+\uref_{l+1,\zeta}^{\ABC}\!\cdot\!\left(\frac{1- \zeta}{2}\right)\!\!{}_{,\zeta}
\dif\zeta&
=0,
\label{E:abc_interface_gradu_1}
\end{align}
where $\alpha_l$ is defined in \eqref{E:alpha_l}.
Noting that the gradient terms are evaluated exactly under $N$-point rule,
one can apply the divergence theorem to obtain:
\begin{align}
\left(\frac{2}{\alpha_l}\right)\gamma\left[
\rint_{-1}^{1}\left(\left(\frac{\alpha_l}{2}\right)^2
\uref_l^{\ABC}-\uref_{l,\zeta\zeta}^{\ABC}\right)
\!\cdot\!
\left(\frac{1+\zeta}{2}\right)
\dif\zeta
+\uref_{l,\zeta}^{\ABC}(+1)
\right]&\notag\\
+\left(\frac{2}{\alpha_{l+1}}\right)\gamma\left[
\rint_{-1}^{1}\left(\left(\frac{\alpha_{l+1}}{2}\right)^2
\uref_{l+1}^{\ABC}-\uref_{l+1,\zeta\zeta}^{\ABC}\right)
\!\cdot\!
\left(\frac{1-\zeta}{2}\right)
\dif\zeta
-\uref_{l+1,\zeta}^{\ABC}(-1)
\right]&
=0,
\label{E:abc_interface_gradu_2}
\end{align}
Using \eqref{E:ul} and \eqref{E:ulzz}, \eqref{E:abc_c=c}, and \eqref{E:ulz},
\eqref{E:abc_interface_gradu_2} becomes:
\begin{align}
\left(\frac{2}{\alpha_l}\right)\gamma
\bigg[
\int_{-1}^{1}
\left(\frac{\alpha_l}{2}\right)^2c^{\ABC}_{l,N-1}\alpha_l^{N-1}\frac{(N+1)!}{(2N)!}P_{N-1}^{(1,1)}(\zeta)
\!\cdot\!
\left(\frac{1+\zeta}{2}\right)
\dif\zeta&
\notag\\
+\rint_{-1}^{1}\left(\frac{\alpha_l}{2}\right)^2c^{\ABC}_{l,N}\alpha_l^{N}\frac{N!}{(2N)!}P_{N}^{(0,0)}(\zeta)
\!\cdot\!
\left(\frac{1+\zeta}{2}\right)
\dif\zeta&
\notag\\
+\left(\frac{\alpha_l}{2}\right)
\sum_{n=0}^{N-1}c^{\ABC}_{l,n+1}\alpha_l^n\frac{(2N- n)!}{(2N)!}P_n^{(N- n,N- n)}(+1)&\bigg]
\notag\\
+\left(\frac{2}{\alpha_{l+1}}\right)\gamma
\bigg[
\int_{-1}^{1}
\left(\frac{\alpha_{l+1}}{2}\right)^2c^{\ABC}_{l+1,N-1}\alpha_{l+1}^{N-1}\frac{(N+ 1)!}{(2N)!}P_{N- 1}^{(1,1)}(\zeta)
\!\cdot\!
\left(\frac{1-\zeta}{2}\right)
\dif\zeta&
\notag\\
+\rint_{-1}^{1}
\left(\frac{\alpha_{l+1}}{2}\right)^2c^{\ABC}_{l+1,N}\alpha_{l+1}^{N}\frac{N!}{(2N)!}P_{N}^{(0,0)}(\zeta)
\!\cdot\!
\left(\frac{1-\zeta}{2}\right)
\dif\zeta&
\notag\\
-\left(\frac{\alpha_{l+1}}{2}\right)
\sum_{n=0}^{N-1}c^{\ABC}_{l+1,n+1}\alpha_{l+1}^n\frac{(2N- n)!}{(2N)!}P_n^{(N- n,N- n)}(-1)&\bigg]=0.
\label{E:abc_interface_gradu_3}%
\end{align}
The second term in each bracket vanishes due to the reduced integration.
Using
(7.391)\textsubscript{4} and
(8.961)\textsubscript{1} of \cite{Gradshtein2014},
we have:
\begin{subequations}
\begin{align}
\int_{-1}^{1}(1+ \zeta)P_{N-1}^{(1,1)}(\zeta)\dif\xi&=\frac{4}{N+1},\\
\int_{-1}^{1}(1- \zeta)P_{N-1}^{(1,1)}(\zeta)\dif\xi&=(- 1)^{N-1}\frac{4}{N+1},
\end{align}
\label{E:P_int}%
\end{subequations}
and, using
(8.960)\textsubscript{1} and
(8.961)\textsubscript{1} of \cite{Gradshtein2014},
we have:
\begin{subequations}
\begin{align}
P_n^{(N-n,N-n)}(+ 1)&={{N}\choose{n}},\\
P_n^{(N-n,N-n)}(- 1)&=(- 1)^n{{N}\choose{n}}.
\end{align}
\label{E:P_at_+-1}%
\end{subequations}
Using \eqref{E:P_int} and \eqref{E:P_at_+-1},
and noting \eqref{E:abc_c=c},
\eqref{E:abc_interface_gradu_3} simplifies to:
\begin{align}
\gamma\left(
S_{\textrm{odd}}^{\ABC}(\alpha_l)\!\cdot\!c_{l,\textrm{even}}^{\ABC}
+ S_{\textrm{even}}^{\ABC}(\alpha_l)\!\cdot\!c_{l,\textrm{odd}}^{\ABC}
\right)
+
\gamma\left(
S_{\textrm{odd}}^{\ABC}(\alpha_{l+1})\!\cdot\!c_{l+1,\textrm{even}}^{\ABC}
- S_{\textrm{even}}^{\ABC}(\alpha_{l+1})\!\cdot\!c_{l+1,\textrm{odd}}^{\ABC}
\right)
=0,
\label{E:abc_interface_gradu_4}
\end{align}
where:
\begin{subequations}
\begin{align}
S_{\textrm{even}}^{\ABC}(\alpha_l)&:=
\sum_{\substack{n=0\\n\textrm{ even}}}^{N}
\frac{\alpha_l^n}{n!}
\frac{N!/(N-n)!}{(2N)!/(2N-n)!},\\
S_{\textrm{odd}}^{\ABC}(\alpha_l)&:=
\sum_{\substack{n=0\\n\textrm{ odd}}}^{N}
\frac{\alpha_l^n}{n!}
\frac{N!/(N-n)!}{(2N)!/(2N-n)!}.
\end{align}
\label{E:abc_S}%
\end{subequations}
Note that:
\begin{subequations}
\begin{align}
S_{\textrm{even}}^{\ABC}(\alpha_l)+
S_{\textrm{odd}}^{\ABC}(\alpha_l)
&={}_1F_1(- N;- 2N;+\alpha_l),\\
S_{\textrm{even}}^{\ABC}(\alpha_l)-
S_{\textrm{odd}}^{\ABC}(\alpha_l)
&={}_1F_1(- N;- 2N;-\alpha_l),
\end{align}
\label{E:FN}%
\end{subequations}
where ${}_1F_1$ represents the confluent hypergeometric function
defined as
(cf. (9.210)\textsubscript{1} of \cite{Gradshtein2014}):
\begin{align}
{}_1F_1(p;q;\beta)=
1
+\frac{p}{q}\frac{\beta^1}{1!}
+\frac{p(p+1)}{q(q+1)}\frac{\beta^2}{2!}
+\frac{p(p+1)(p+2)}{q(q+1)(q+2)}\frac{\beta^3}{3!}
+\cdots.
\end{align}
Crucial in this work is the relation
of the confluent hypergeometric functions
to the Pad\'{e} approximants; specifically:
\begin{align}
[N/N]_{\exp{(z)}}=\frac{{}_1F_1(-N;-2N;+z)}{{}_1F_1(-N;-2N;-z)},
\label{E:nnpade}
\end{align}
where $[N/N]_{\exp{(z)}}$ is the Pad\'{e} approximant of $\exp{(z)}$
of order $[N/N]$;
see, e.g., \cite{BirkhoffVarga1965,SaffVarga1975,Ainsworth2004}.
In the following we use a shorthand notation,
${}_1F_1^N(z)$, for ${}_1F_1(-N;-2N;z)$.

On the other hand, in the $\zeta$-coordinate system,
\eqref{E:abc_interface_u} is written as:
\begin{align}
\uref^{\ABC}_l(+1)-\uref^{\ABC}_{l+1}(-1)=0.
\label{E:abc_interface_u_1}
\end{align}
Using \eqref{E:ul}, \eqref{E:P_at_+-1}, \eqref{E:abc_c=c},
and \eqref{E:abc_S},
one can write $\uref_l^{\ABC}(+1)$ and $\uref_{l}^{\ABC}(-1)$ as:
\begin{subequations}
\begin{align}
\uref_l^{\ABC}(+1)&=
S_{\textrm{even}}^{\ABC}(\alpha_l)\!\cdot\!c_{l,\textrm{even}}^{\ABC}
+
S_{\textrm{odd}}^{\ABC}(\alpha_l)\!\cdot\!c_{l,\textrm{odd}}^{\ABC},
\label{E:abc_interface_ul_left}\\
\uref_l^{\ABC}(-1)&=
S_{\textrm{even}}^{\ABC}(\alpha_l)\!\cdot\!c_{l,\textrm{even}}^{\ABC}
-
S_{\textrm{odd}}^{\ABC}(\alpha_l)\!\cdot\!c_{l,\textrm{odd}}^{\ABC}.
\label{E:abc_interface_ul_right}
\end{align}
\label{E:abc_interface_ul}%
\end{subequations}
Using \eqref{E:abc_interface_ul},
\eqref{E:abc_interface_u_1} is then written as:
\begin{align}
\left(
S_{\textrm{even}}^{\ABC}(\alpha_l)\!\cdot\!c_{l,\textrm{even}}^{\ABC}
+ S_{\textrm{odd}}^{\ABC}(\alpha_l)\!\cdot\!c_{l,\textrm{odd}}^{\ABC}
\right)
-
\left(
S_{\textrm{even}}^{\ABC}(\alpha_{l+1})\!\cdot\!c_{l+1,\textrm{even}}^{\ABC}
- S_{\textrm{odd}}^{\ABC}(\alpha_{l+1})\!\cdot\!c_{l+1,\textrm{odd}}^{\ABC}
\right)
=0,
\label{E:abc_interface_u_2}
\end{align}

Using
\eqref{E:abc_interface_u_2} and
\eqref{E:abc_interface_gradu_4},
one can now set up the following system:
\begin{align}
\begin{pmatrix}
S_{\textrm{even}}^{\ABC}(\alpha_l) & S_{\textrm{odd}}^{\ABC}(\alpha_l)\\
S_{\textrm{odd}}^{\ABC}(\alpha_l) & S_{\textrm{even}}^{\ABC}(\alpha_l)
\end{pmatrix}
\hspace{-5pt}
\begin{pmatrix}
c_{l,\textrm{even}}^{\ABC}\\
c_{l,\textrm{odd}}^{\ABC}
\end{pmatrix}
-
\begin{pmatrix}
\phantom{-\hspace{-2pt}}S_{\textrm{even}}^{\ABC}(\alpha_{l+1}) &
\hspace{-10pt}-\hspace{-2pt}S_{\textrm{odd}}^{\ABC}(\alpha_{l+1})\\
-\hspace{-2pt}S_{\textrm{odd}}^{\ABC}(\alpha_{l+1}) &
\hspace{-10pt}\phantom{-\hspace{-2pt}}S_{\textrm{even}}^{\ABC}(\alpha_{l+1})
\end{pmatrix}
\hspace{-5pt}
\begin{pmatrix}
c_{l+1,\textrm{even}}^{\ABC}\\
c_{l+1,\textrm{odd}}^{\ABC}
\end{pmatrix}
=
\begin{pmatrix}
0\\
0
\end{pmatrix}
.
\label{E:abc_interface}
\end{align}
Performing eigendecompositions of the matrices appearing
in \eqref{E:abc_interface}, and using \eqref{E:FN}, we have:
\begin{subequations}
\begin{align}
\begin{pmatrix}
S_{\textrm{even}}^{\ABC}(\alpha_l) & S_{\textrm{odd}}^{\ABC}(\alpha_l)\\
S_{\textrm{odd}}^{\ABC}(\alpha_l) & S_{\textrm{even}}^{\ABC}(\alpha_l)
\end{pmatrix}
&=
\Phi
\begin{pmatrix}
{}_1F_1^N(+\alpha_l) & \hspace{-15pt}0\\
0 & \hspace{-15pt}{}_1F_1^N(-\alpha_l)
\end{pmatrix}
\Phi^\textrm{T},\\
\begin{pmatrix}
\phantom{-\hspace{-2pt}}S_{\textrm{even}}^{\ABC}(\alpha_{l+1}) &
\hspace{-10pt}-\hspace{-2pt}S_{\textrm{odd}}^{\ABC}(\alpha_{l+1})\\
-\hspace{-2pt}S_{\textrm{odd}}^{\ABC}(\alpha_{l+1}) &
\hspace{-10pt}\phantom{-\hspace{-2pt}}S_{\textrm{even}}^{\ABC}(\alpha_{l+1})
\end{pmatrix}
&=
\Phi
\begin{pmatrix}
{}_1F_1^N(-\alpha_{l+1}) & \hspace{-15pt}0\\
0 & \hspace{-15pt}{}_1F_1^N(+\alpha_{l+1})
\end{pmatrix}
\Phi^\textrm{T},
\end{align}
\label{E:abc_eigen}%
\end{subequations}
where $\Phi$ is defined in \eqref{E:Phi}.
Using \eqref{E:abc_eigen} and \eqref{E:Phi},
\eqref{E:abc_interface} reduces to:
\begin{align}
\begin{pmatrix}
{}_1F_1^N(+\alpha_l) & \hspace{-15pt}0\\
0 & \hspace{-15pt}{}_1F_1^N(-\alpha_l)
\end{pmatrix}
\hspace{-5pt}
\begin{pmatrix}
c_{l,-}^{\ABC}\\
c_{l,+}^{\ABC}
\end{pmatrix}
=
\begin{pmatrix}
{}_1F_1^N(-\alpha_{l+1}) & \hspace{-15pt}0\\
0 & \hspace{-15pt}{}_1F_1^N(+\alpha_{l+1})
\end{pmatrix}
\hspace{-5pt}
\begin{pmatrix}
c_{l+1,-}^{\ABC}\\
c_{l+1,+}^{\ABC}
\end{pmatrix}
,
\label{E:abc_recursion}
\end{align}
where:
\begin{align}
\begin{pmatrix}
c_{l,-}^{\ABC}\\
c_{l,+}^{\ABC}
\end{pmatrix}
=
\Phi^\textrm{T}
\begin{pmatrix}
c_{l,\textrm{even}}^{\ABC}\\
c_{l,\textrm{odd}}^{\ABC}
\end{pmatrix}
.
\label{E:abc_cl+-}
\end{align}
Using \eqref{E:abc_recursion} for $l\in\{0,\cdots,L-1\}$, one obtains:
\begin{subequations}
\begin{align}
\left\{
\begin{array}{rl}
c^{\ABC}_{l,-}=0,
\phantom{c_{L,-}^{\ABC}}
&
\phantom{l_{\textrm{max}}}0\leq l<l_{\textrm{max}},\\
\prod\limits_{l'=l}^{L-1}\left(\frac{{}_1F_1^N(+\alpha_{l'})}{{}_1F_1^N(-\alpha_{l'+1})}\right)
\!\cdot\!
c_{l,-}^{\ABC}
=
c_{L,-}^{\ABC},
\phantom{0}
&
\phantom{0}l_{\textrm{max}}\leq l<L,
\end{array}
\hspace{16pt}
\right.\\
\intertext{and:}
\left\{
\begin{array}{rl}
\prod\limits_{l'=l}^{l_{\textrm{min}}-1}\left(\frac{{}_1F_1^N(-\alpha_{l'})}{{}_1F_1^N(+\alpha_{l'+1})}\right)
\!\cdot\!
c_{l,+}^{\ABC}
=
c_{l_{\textrm{min}},+}^{\ABC},
\phantom{0}
&
\phantom{l_{\textrm{min}}}0\leq l<l_{\textrm{min}},\\
c^{\ABC}_{l,+}=0,
\phantom{c_{l_{\textrm{min}},+}^{\ABC}}
&
\phantom{0}l_{\textrm{min}}<l\leq L,
\end{array}
\right.
\end{align}
\label{E:abc_c}%
\end{subequations}
where:
\begin{subequations}
\begin{align}
l_{\textrm{max}}&:=\phantom{\max}\hspace{-20pt}
\min
\{
l\in\{0,\cdots,L\}|{}_1F_1^N(-\alpha_{l'})\neq 0,\forall l'>l
\},\\
l_{\textrm{min}}&:=\phantom{\min}\hspace{-20pt}
\max
\{
l\in\{0,\cdots,L\}|{}_1F_1^N(-\alpha_{l'})\neq 0,\forall l'<l
\},
\end{align}
\label{E:abc_l}%
\end{subequations}
leaving only two undetermined coefficients,
$c^{\ABC}_{L,-}$ and $c^{\ABC}_{0,+}$
Note that $l_{\textrm{max}}=0$ and $l_{\textrm{min}}=L$
almost everywhere in $\real^{d-1}\ni(k_2,\cdots,k_d)$ for a given value of $s$.

$c_{0,-}^{\ABC}$ and $c_{0,+}^{\ABC}$ can be viewed as representing
the incoming/growing and the outgoing/decaying components of
the \emph{discrete wave} in $x<0$
in the neighbourhood of the boundary.
Here, we suppose that there was an additional layer $[x_{-2},x_{-1}]$ with
$-\delta<x_{-2}<x_{-1}<0$ and
define $h_{-1}:=x_{-1}-x_{-2}=h_0$ and
$\gamma_{-1}:=1=\gamma_0$.
If we extended the above analysis to include $[x_{-2},x_{-1}]$,
also using the reduced integration in $[x_{-2},x_{-1}]$,
$u_{-1}^{\ABC}(x)$ on $[x_{-2},x_{-1}]$
would be represented in terms of
$c_{-1,-}^{\ABC}$ and $c_{-1,+}^{\ABC}$.
Using \eqref{E:abc_recursion} with $l=-1$, and noting \eqref{E:alpha_l},
$c_{-1,-}^{\ABC}$ and $c_{-1,+}^{\ABC}$ would be
related to $c_{0,-}^{\ABC}$ and $c_{0,+}^{\ABC}$ as:
\begin{align}
\begin{pmatrix}
{}_1F_1^N(+\gamma h_0) & \hspace{-15pt}0\\
0 & \hspace{-15pt}{}_1F_1^N(-\gamma h_0)
\end{pmatrix}
\hspace{-5pt}
\begin{pmatrix}
c_{-1,-}^{\ABC}\\
c_{-1,+}^{\ABC}
\end{pmatrix}
=
\begin{pmatrix}
{}_1F_1^N(-\gamma h_0) & \hspace{-15pt}0\\
0 & \hspace{-15pt}{}_1F_1^N(+\gamma h_0)
\end{pmatrix}
\hspace{-5pt}
\begin{pmatrix}
c_{0,-}^{\ABC}\\
c_{0,+}^{\ABC}
\end{pmatrix}
.
\end{align}
By \eqref{E:nnpade},
${}_1F_1^N(+\gamma h_0)/{}_1F_1^N(-\gamma h_0)$ and
${}_1F_1^N(-\gamma h_0)/{}_1F_1^N(+\gamma h_0)$ are, respectively, the
Pad\'{e} approximants of $\exp{(+\gamma h_0)}$ and $\exp{(-\gamma h_0)}$ of order $[N/N]$,
which implies the above statement.

\subsection{Termination condition}\label{SS:abc_end}
Using \eqref{E:abc_interface_ul_left} with $l=L$,
\eqref{E:abc_end_u} can be written as:
\begin{align}
S_{\textrm{even}}^{\ABC}(\alpha_L)\!\cdot\!c_{L,\textrm{even}}^{\ABC}
+ S_{\textrm{odd}}^{\ABC}(\alpha_L)\!\cdot\!c_{L,\textrm{odd}}^{\ABC}
=0.
\label{E:abc_end_u_1}
\end{align}
Rewriting in terms of $c_{L,-}^{\ABC}$ and $c_{L,+}^{\ABC}$
using \eqref{E:abc_cl+-}, and noting \eqref{E:FN},
we have:
\begin{align}
\frac{1}{\sqrt{2}}
\begin{pmatrix}
{}_1F_1^N(+\alpha_L)\hspace{5pt}{}_1F_1^N(-\alpha_L)
\end{pmatrix}
\hspace{-3pt}
\begin{pmatrix}
c_{L,-}^{\ABC}\\
c_{L,+}^{\ABC}
\end{pmatrix}
=0.
\label{E:abc_end_u_2}
\end{align}
Using \eqref{E:abc_c} in conjunction with \eqref{E:abc_end_u_2},
one is left with one undetermined coefficient, $c^{\ABC}_{0,+}$, as desired.
Specifically, one obtains:
\begin{align}
c_{0,-}^{\ABC}=
-
\left(
\frac{{}_1F_1^N(-\alpha_0)}{{}_1F_1^N(+\alpha_0)}
\right)
\prod_{l=1}^{L}
\left(
\frac{{}_1F_1^N(-\alpha_l)}{{}_1F_1^N(+\alpha_l)}
\right)^2
\!\cdot\!c_{0,+}^{\ABC}
,
\label{E:abc_c0-/c0+}
\end{align}
which can be written, using \eqref{E:nnpade}, as:
\begin{align}
c_{0,-}^{\ABC}=
-
\left(
[N/N]_{\exp{(-\alpha_0)}}
\right)
\prod_{l=1}^{L}
\left(
[N/N]_{\exp{(-\alpha_l)}}
\right)^2
\!\cdot\!c_{0,+}^{\ABC}
,
\label{E:abc_reflection}
\end{align}
which compares with \eqref{E:const_reflection}.

\begin{table}
\begin{tabular}{|c|l|}
\hline
$N$ & \hspace{15pt}Zeros of $[N/N]_{\exp{(-z)}}$ \\
\hline
\multirow{1}{*}{1} & $2.00000000$\\ 
\hline
\multirow{2}{*}{2} & $3.00000000-1.73205081\iu$\\ 
                   & $3.00000000+1.73205081\iu$\\ 
\hline
\multirow{3}{*}{3} & $4.64437071$\\ 
                   & $3.67781465-3.50876192\iu$\\ 
                   & $3.67781465+3.50876192\iu$\\ 
\hline
\multirow{4}{*}{4} & $4.20757879-5.31483608\iu$ \\ 
                   & $4.20757879+5.31483608\iu$ \\ 
                   & $5.79242121-1.73446826\iu$ \\ 
                   & $5.79242121+1.73446826\iu$\\ 
\hline
\end{tabular}
\caption{Zeros of $[N/N]_{\exp{(-z)}}$ for $N\in\{1,2,3,4\}$.}
\label{Ta:pade}
\end{table}

It is known that all zeros of $[N/N]_{\exp{(-z)}}$,
denoted as $\{z_{n}\}_{n\in\{1,\cdots,N\}}$,
come in conjugate pairs if not real and
lie in the open right half complex plane, and,
due to the symmetry between the denominator and the numerator in \eqref{E:nnpade},
all poles of $[N/N]_{\exp{(-z)}}$ are given as
$\{-z_{n}\}_{n\in\{1,\cdots,N\}}$;
see \cite{BirkhoffVarga1965,SaffVarga1975} for details, and see
Table \ref{Ta:pade} for the zeros of $[N/N]_{\exp{(-z)}}$ for $N\in\{1,2,3,4\}$.
Recalling that $\Real{\alpha_l}>0$ for $l\in\{1,\cdots,L\}$,
one can readily show that $\abs{[N/N]_{\exp{(-\alpha_l)}}}<1$
for each $l\in\{1,\cdots,L\}$, and
the absolute value of the reflection coefficient appearing in \eqref{E:abc_reflection}
decreases as $O(R^{2LN})$ for some $R<1$, as $L$ and $N$ increase.
Furthermore,
if any $\alpha_l$, $l\in\{1,\cdots,L\}$, coincides with a zero of $[N/N]_{\exp{(-z)}}$,
the outgoing/decaying wave is annihilated at the boundary,
causing no spurious reflection.

We now study the ability of the proposed ABCs to approximate
the Sommerfeld radiation condition \eqref{E:helmholtz_sommerfeld}.
This analysis can be performed without assuming that the $N$-point
reduced integration is also used in the physical domain, $x<0$;
i.e., we do not use \eqref{E:abc_domain} with $l=0$ and \eqref{E:abc_interface_gradu} with $l=0$.

Using \eqref{E:abc_domain} with $l\in\{1,\cdots,L\}$ and
\eqref{E:abc_interface_gradu} with $l\in\{1,\cdots,L-1\}$,
we obtain the following identity that is similar to \eqref{E:abc_c0-/c0+}:
\begin{align}
c_{1,-}^{\ABC}=
-
\left(
\frac{{}_1F_1^N(-\alpha_1)}{{}_1F_1^N(+\alpha_1)}
\right)
\prod_{l=2}^{L}
\left(
\frac{{}_1F_1^N(-\alpha_l)}{{}_1F_1^N(+\alpha_l)}
\right)^2
\!\cdot\!c_{1,+}^{\ABC}
.
\label{E:abc_c1-/c1+}
\end{align}
The approximation to $u_{0,x}^{\ABC}(0)$ provided by the proposed ABCs,
which is weakly enforced to the problem solved in the physical domain $x<0$,
is denoted by $[u_{0,x}^{\ABC}(0)]_{\weak}$, and is given as:
\begin{align}
[{u_{,x}^{\ABC}}(0)]_{\weak}=-
\rint_{x_{0}}^{x_1}\frac{\gamma^2}{\gamma_1}u_1^{\ABC}w_0^{\HAT}
+\gamma_1u_{1,x}^{\ABC}w_{0,x}^{\HAT}\dif x.
\label{E:abc_approx_ux}
\end{align}
Using \eqref{E:abc_interface_ul_right} with $l=1$ along with
\eqref{E:abc_interface_u} with $l=0$ and
using the second half of \eqref{E:abc_interface_gradu_4} with $l=0$, we obtain:
\begin{subequations}
\begin{align}
u^{\ABC}_0(0)&=
S_{\textrm{even}}^{\ABC}(\alpha_1)\!\cdot\!c_{1,\textrm{even}}^{\ABC}-
S_{\textrm{odd}}^{\ABC}(\alpha_1)\!\cdot\!c_{1,\textrm{odd}}^{\ABC},\\
[u_{0,x}^{\ABC}(0)]_{\weak}&=
-\gamma\left(
S_{\textrm{odd}}^{\ABC}(\alpha_1)\!\cdot\!c_{1,\textrm{even}}^{\ABC}
- S_{\textrm{even}}^{\ABC}(\alpha_1)\!\cdot\!c_{1,\textrm{odd}}^{\ABC}
\right)
.
\end{align}
\label{E:abc_approx_1}%
\end{subequations}
Using \eqref{E:abc_cl+-}, and noting \eqref{E:FN}, one can rewrite
\eqref{E:abc_approx_1} in terms of
$c_{1,-}^{\ABC}$ and $c_{1,+}^{\ABC}$ as:
\begin{subequations}
\begin{align}
u^{\ABC}_0(0)&=
\frac{1}{\sqrt{2}}
\begin{pmatrix}
{}_1F_1^N(-\alpha_1)\hspace{5pt}{}_1F_1^N(+\alpha_1)
\end{pmatrix}
\hspace{-3pt}
\begin{pmatrix}
c_{1,-}^{\ABC}\\
c_{1,+}^{\ABC}
\end{pmatrix}
.\\
[u_{0,x}^{\ABC}(0)]_{\weak}&=
-\frac{1}{\sqrt{2}}\gamma
\begin{pmatrix}
-\hspace{-2pt}{}_1F_1^N(-\alpha_1)\hspace{5pt}{}_1F_1^N(+\alpha_1)
\end{pmatrix}
\hspace{-3pt}
\begin{pmatrix}
c_{1,-}^{\ABC}\\
c_{1,+}^{\ABC}
\end{pmatrix}
.
\end{align}
\label{E:abc_approx_2}%
\end{subequations}
Using \eqref{E:abc_approx_2} in conjunction with \eqref{E:abc_c1-/c1+},
we obtain:
\begin{align}
[{u_{0,x}^{\ABC}}(0)]_{\weak}
+\gamma
\left(
\frac{1+
\prod\limits_{l=1}^{L}
\left(
\frac{{}_1F_1^N(-\alpha_l)}{{}_1F_1^N(+\alpha_l)}
\right)^2
}{1-
\prod\limits_{l=1}^{L}
\left(
\frac{{}_1F_1^N(-\alpha_l)}{{}_1F_1^N(+\alpha_l)}
\right)^2
}
\right)
\hspace{-2pt}
{u^{\ABC}_0}(0)=0,
\end{align}
which can be rewritten noting \eqref{E:nnpade} as:
\begin{align}
[{u_{0,x}^{\ABC}}(0)]_{\weak}
+\gamma
\left(
\frac{1+
\prod\limits_{l=1}^{L}
\left(
[N/N]_{\exp{(-\alpha_l)}}
\right)^2
}{1-
\prod\limits_{l=1}^{L}
\left(
[N/N]_{\exp{(-\alpha_l)}}
\right)^2
}
\right)
\hspace{-2pt}
u^{\ABC}_0(0)=0,
\label{E:abc_approx}
\end{align}
which compares with \eqref{E:const_sommerfeld}.
As noted in the above, the product terms in \eqref{E:abc_approx}
decrease as $O(R^{2LN})$ for some $R<1$ as $L$ and $N$ increase, and
\eqref{E:abc_approx} approaches to
the exact Sommerfeld radiation condition \eqref{E:helmholtz_sommerfeld}.
The approximation obtained here, i.e., \eqref{E:abc_approx},
generalises that obtained in \cite{GuddatiLim2006} for the PMDLs.

\subsection{Transformed Sommerfeld radiation condition}
The analysis performed in
Sec.~\ref{SS:abc_domain}, Sec.~\ref{SS:abc_interface}, and Sec.~\ref{SS:abc_end}
can readily be modified for a problem in which
we apply the transformed Sommerfeld radiation condition,
$\gamma_Lu_{L,x}^{\ABC}(x_L)+\gamma u_L^{\ABC}(x_L)=0$, at $x=x_L$
weakly as the Neumann boundary condition,
in place of the termination condition \eqref{E:abc_end_u}.
Specifically, we enforce the following condition instead of \eqref{E:abc_end_u}:
\begin{align}
\rint_{x_{L-1}}^{x_L}\frac{\gamma^2}{\gamma_L}u_L^{\ABC}w_L^{\HAT}
+\gamma_Lu_{L,x}^{\ABC}w_{L,x}^{\HAT}\dif x
-
(-\gamma u_L^{\ABC}(x_L)w_L^{\HAT}(x_L))
=0,
\label{E:abc_end_ux}
\end{align}
Using the first half of \eqref{E:abc_interface_gradu_4} with $l=L$ and
\eqref{E:abc_interface_ul_left} with $l=L$,
one can rewrite \eqref{E:abc_end_ux} as:
\begin{align}
\gamma\left(
S_{\textrm{odd}}^{\ABC}(\alpha_L)\!\cdot\!c_{L,\textrm{even}}^{\ABC}
+ S_{\textrm{even}}^{\ABC}(\alpha_L)\!\cdot\!c_{L,\textrm{odd}}^{\ABC}
\right)
+
\gamma\left(
S_{\textrm{even}}^{\ABC}(\alpha_L)\!\cdot\!c_{L,\textrm{even}}^{\ABC}+
S_{\textrm{odd}}^{\ABC}(\alpha_L)\!\cdot\!c_{L,\textrm{odd}}^{\ABC}
\right)
=0,
\label{E:abc_end_ux_1}
\end{align}
Using \eqref{E:abc_cl+-}, and noting \eqref{E:FN},
one can rewrite \eqref{E:abc_end_ux_1}
in terms of $c_{L,-}^{\ABC}$ and $c_{L,+}^{\ABC}$ as:
\begin{align}
\sqrt{2}\gamma
\begin{pmatrix}
{}_1F_1^N(+\alpha_L)\hspace{5pt}0
\end{pmatrix}
\hspace{-3pt}
\begin{pmatrix}
c_{L,-}^{\ABC}\\
c_{L,+}^{\ABC}
\end{pmatrix}
=0.
\label{E:abc_end_ux_2}
\end{align}
Using \eqref{E:abc_c} in conjunction with \eqref{E:abc_end_ux_2},
instead of with \eqref{E:abc_end_u_2},
we obtain:
\begin{align}
c_{0,-}^{\ABC}=0\!\cdot\!c_{0,+}^{\ABC},
\label{E:abc_reflection_ux}
\end{align}
which compares with \eqref{E:const_sommerfeld_transformed_reflection}, and
implies that \eqref{E:abc} with \eqref{E:abc_end_u}
replaced with \eqref{E:abc_end_ux}
only admits the outgoing/decaying wave solution.

Following the same lines as the discussion in Sec.~\ref{SS:abc_end},
one can show,
without assuming that the $N$-point reduced integration is also used
in the physical domain, that:
\begin{align}
[{u_{0,x}^{\ABC}}(0)]_{\weak}
+\gamma
u^{\ABC}_0(0)=0,
\label{E:abc_approx_ux}
\end{align}
which compares with \eqref{E:const_sommerfeld_transformed_sommerfeld},
and implies that the proposed ABCs with
the transformed Sommerfeld radiation condition at $x=x_L$
would provide the exact Sommerfeld radiation condition
\eqref{E:helmholtz_sommerfeld} at $x=0$
for the problem solved in the physical domain.
This property has been discovered for the special case of the PMDLs
in \cite{GuddatiLim2006},
though we used slightly different arguments.

\section{Numerical examples}\label{S:example}
We solved one- and three-dimensional example problems numerically
using Firedrake \cite{FiredrakeUserManual,Rathgeber2016}.
Firedrake is a high-level, high-productivity system for
the specification and solution of partial differential equations
using finite element methods.
Firedrake uses PETSc \cite{petsc-web-page,petsc-user-ref,petsc-efficient}
for the mesh representations and the linear/nonlinear solvers.
Firedrake has been used for implementing PMLs in \cite{Kirby2021}.

\begin{figure}[ht]
\centering
\includegraphics[scale=1]{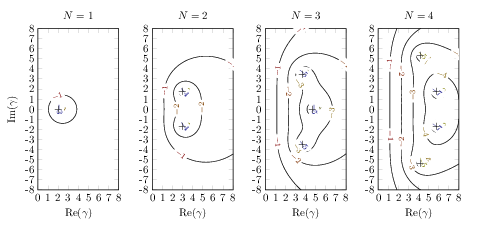}
\caption{
Reflection coefficient
computed for the one-dimensional Helmholtz problem
with $\gamma\in\{z\in\complex|0<\Real{z}<8,-8<\Imag{z}<+8\}$
in $(-1,0)$
with the ABC of type $(1,N)$ with $\gamma_1=1$
applied in $(0,1)$,
where $N\in\{1,2,3,4\}$.
A $Q_N$ Lagrange finite element with $N$-point reduced integration
was used in $(-1,0)$ as well as in $(0,1)$.
For each $N\in\{1,2,3,4\}$,
the reflection coefficient virtually vanished where $\gamma/\gamma_1$
coincided with the zeros of $[N/N]_{\exp{(-z)}}$;
those points are shown with the $+$ symbols.
}
\label{Fi:one_dim}
\end{figure}

\begin{figure}[ht]
\centering
\includegraphics[scale=1]{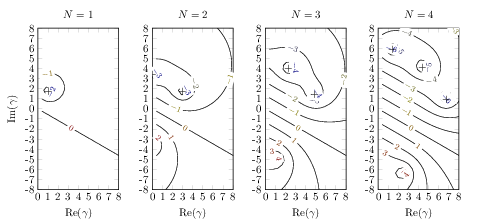}
\caption{
Reflection coefficient
computed for the one-dimensional Helmholtz problem
with $\gamma\in\{z\in\complex|0<\Real{z}<8,-8<\Imag{z}<+8\}$
in $(-1,0)$
with the ABC of type $(1,N)$ with $\gamma_1=1/2+\sqrt{3}/2\!\cdot\!\iu$
applied in $(0,1)$,
where $N\in\{1,2,3,4\}$.
A $Q_N$ Lagrange finite element with $N$-point reduced integration
was used in $(-1,0)$ as well as in $(0,1)$.
For each $N\in\{1,2,3,4\}$,
the reflection coefficient virtually vanished where $\gamma/\gamma_1$
coincided with the zeros of $[N/N]_{\exp{(-z)}}$;
those points are shown with the $+$ symbols.
}
\label{Fi:one_dim_rotated}
\end{figure}

\subsection{One-dimensional examples}
In this section we numerically verify the formula for the reflection coefficient
\eqref{E:abc_reflection} for ABCs of type $(1,N)$.

We consider a one-dimensional domain, $D^{\oneD}:=(-1,1)$, in which
the physical domain and the artificial domain are defined to be
$(-1,0)$ and $(0,1)$, respectively.
We consider a mesh $\mathcal{T}^{\oneD}$ of $D^{\oneD}$
defined as $\mathcal{T}^{\oneD}:=\{[-1,0], [0,1]\}$;
i.e., we have one absorbing layer, $[0,1]$.
We then define function spaces on $\mathcal{T}^{\oneD}$,
$V^{\oneD}$ and $V^{\oneD}_0$, as:
\begin{subequations}
\begin{align}
V^{\oneD}:&=\{w\in H^1(D^{\oneD})|w_{|K}\in\mathbb{P}_N,\forall K\in\mathcal{T}^{\oneD}\},\\
V^{\oneD}_0:&=\{w\in V^{\oneD}|w(\pm 1)=0\},
\end{align}
\end{subequations}
where $\mathbb{P}_N$ is the space of univariate polynomials of degree $N$.
We use the $Q_N$, or $P_N$, Lagrange finite element
with the Gauss-Lobatto nodes.
We then consider the following problem:
\begin{subequations}
\begin{align}
\intertext{Find $u^{\oneD}\in V^{\oneD}$ such that}
\rint_{D^{\oneD}}\frac{\gamma^2}{g^{\oneD}}u^{\oneD}\ol{w^{\oneD}}
+g^{\oneD}u_{,x}^{\oneD}\ol{w_{,x}^{\oneD}}\dif x&=0,
\quad\forall w^{\oneD}\in V^{\oneD}_0,
\label{E:example_1d_domain}\\
u^{\oneD}(+1)&=0,\label{E:example_1d_end}\\
u^{\oneD}(-1)&=1,\label{E:example_1d_dirichlet}
\end{align}
\label{E:example_id}%
\end{subequations}
where:
\begin{align}
g^{\oneD}(x):=\left\{
\begin{array}{cl}
\gamma_0, &x\leq 0,\\
\gamma_1, &\textrm{otherwise},
\end{array}
\right.
\end{align}
where $\gamma_0=1$ and $\gamma_1\in\complex$.
Note that we use $N$-point reduced integration
in the physical domain $(-1,0)$ as well as in the artificial domain $(0,1)$.
\eqref{E:example_1d_end} is the termination condition for the proposed ABC and
\eqref{E:example_1d_dirichlet} is the standard Dirichlet boundary condition.

We solved \eqref{E:example_id}
for $\gamma_1\in\{1,1/2+\sqrt{3}/2\!\cdot\!\iu\}$
for $N\in\{1,2,3,4\}$
for $\gamma\in\{z\in\complex|0<\Real{z}<8,-8<\Imag{z}<+8\}$
using Firedrake.
We used MUMPS parallel sparse direct solver
\cite{mumps1,mumps2} via PETSc.
For each case, to compute the the reflection coefficient
for the discrete wave solution
defined in \eqref{E:abc_reflection} numerically,
we decomposed the computed solution into
the incoming/growing discrete wave component and
the outgoing/decaying discrete wave component
directly inspecting the finite element local matrix
associated with $K=[-1,0]\in\mathcal{T}^{\oneD}$
required to assemble \eqref{E:example_1d_domain};
see \cite{Koyama2008,Sagiyama2013} for details.
Fig.~\ref{Fi:one_dim} and Fig.~\ref{Fi:one_dim_rotated} show the
contour plots of the computed reflection coefficients
for $N\in\{1,2,3,4\}$,
for $\gamma_1=1$ and $\gamma_1=1/2+\sqrt{3}/2\!\cdot\!\iu$, respectively.
One can observe that these contour plots agree with
the formula \eqref{E:abc_reflection}.
Specifically,
for each $\gamma_1$, for each $N$,
we verified
that the reflection coefficient vanished at which
$\gamma/\gamma_1$ coincided with the zeros of $[N/N]_{\exp{(-z)}}$;
Table \ref{Ta:pade} lists the zeros of $[N/N]_{\exp{(-z)}}$
for $N\in\{1,2,3,4\}$.

% Case 0
% N = 1

% N = 2

% N = 3

% N = 4

\begin{figure}[ht]
\centering
\includegraphics[scale=.9]{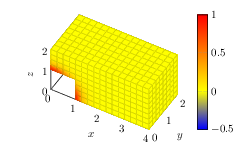}
\caption{Exact solution to the example Helmholtz problem with $s=+4.00+0.25\iu$.}
\label{Fi:example_3d_exact_0}

\vs

\includegraphics[scale=.9]{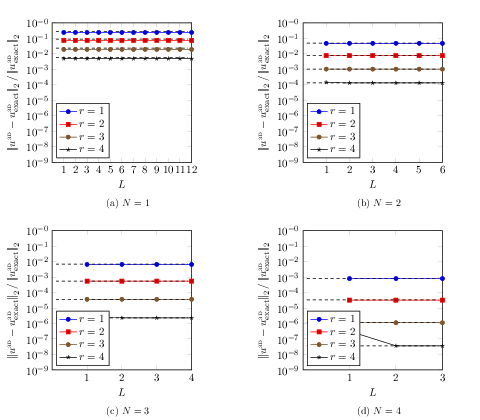}
\caption{
Relative $L^2$-norm errors from ABCs of type $(L,N)$
for $s=+4.00+0.25\iu$:
(a) $N=1$,
(b) $N=2$,
(c) $N=3$, and
(d) $N=4$.
$r\in\{1,2,3,4\}$ represents the refinement level.
$L$ was chosen for each $N$
so that $(LN)_\textrm{max}\leq 12$.
Estimated discretisation errors are also shown with black dashed lines.
}
\label{Fi:convergence_case_0}
\end{figure}

% Case 1
% N = 1

% N = 2

% N = 3

% N = 4

\begin{figure}[ht]
\centering
\includegraphics[scale=.9]{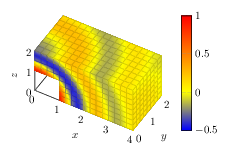}
\caption{Exact solution to the example Helmholtz problem with $s=+0.25+4.00\iu$.}
\label{Fi:example_3d_exact_1}

\vs

\includegraphics[scale=.9]{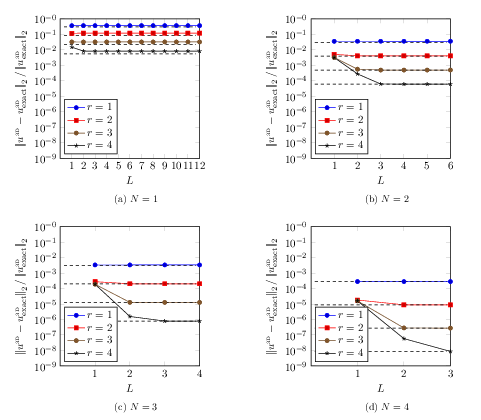}
\caption{
Relative $L^2$-norm errors from ABCs of type $(L,N)$
for $s=+0.25+4.00\iu$:
(a) $N=1$,
(b) $N=2$,
(c) $N=3$, and
(d) $N=4$.
$r\in\{1,2,3,4\}$ represents the refinement level.
$L$ was chosen for each $N$
so that $(LN)_\textrm{max}\leq 12$.
Estimated discretisation errors are also shown with black dashed lines.
}
\label{Fi:convergence_case_1}
\end{figure}

% Case 2
% N = 1

% N = 2

% N = 3

% N = 4

\begin{figure}[ht]
\centering
\includegraphics[scale=.9]{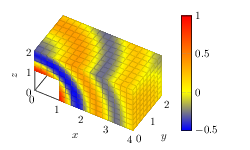}
\caption{Exact solution to the example Helmholtz problem with $s=+4.00\iu$.}
\label{Fi:example_3d_exact_2}

\vs

\includegraphics[scale=.9]{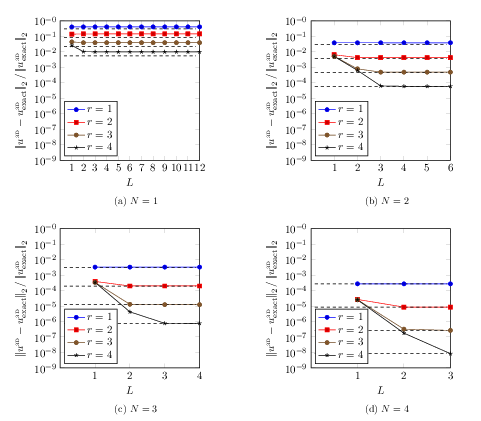}
\caption{
Relative $L^2$-norm errors from ABCs of type $(L,N)$
for $s=+4.00\iu$:
(a) $N=1$,
(b) $N=2$,
(c) $N=3$, and
(d) $N=4$.
$r\in\{1,2,3,4\}$ represents the refinement level.
$L$ was chosen for each $N$
so that $(LN)_\textrm{max}\leq 12$.
Estimated discretisation errors are also shown with black dashed lines.
}
\label{Fi:convergence_case_2}
\end{figure}

\subsection{Three-dimensional examples}
In this section we follow the standard approach
\cite{Berenger1994,GuddatiLim2006}
to extend the one-dimensional formulation
of the proposed ABC introduced in Sec.~\ref{S:abc}
to three dimensions.

We consider a three-dimensional ($d=3$) domain $D^{\threeD}$ defined as
$D^{\threeD}:=\{\bs{x}\in\real^3|\bs{x}\in\prod_{i\in\{1,2,3\}}(0,B_i+Lh)\textrm{ and }\bs{x}\not\in[0,1]^3\}$,
representing a box with a cubic hole,
where $(B_1,B_2,B_3)=(4,2,2)$ and
$h=2^{-\nref}$ with $\nref\in\mathbb{N}$ being the refinement level.
The physical domain is given as
$\{\bs{x}\in\real^3|\bs{x}\in\prod_{i\in\{1,2,3\}}(0,B_i)\textrm{ and }\bs{x}\not\in[0,1]^3\}$, and
we are to use an ABC of type $(L,N)$ in the artificial domain
in each direction.
We let
$\partial D^{\threeD}_0:=\{\bs{x}\in\partial D^{\threeD}|x_i=B_i+Lh\textrm{ for some }i\in\{0,1,2\}\}$,
$\partial D^{\threeD}_1:=\{\bs{x}\in\partial D^{\threeD}|\bs{x}\in[0,1]^3\}$, and
$\partial D^{\threeD}_2:=\{\bs{x}\in\partial D^{\threeD}|x_i=0\textrm{ for some }i\in\{0,1,2\}\}$.
We consider a mesh $\mathcal{T}^{\threeD}$ of $D^{\threeD}$
composed of cuboids of edge lengths $h$.
We then define function spaces on $\mathcal{T}^{\threeD}$,
$V^{\threeD}$ and $V^{\threeD}_0$, as:
\begin{subequations}
\begin{align}
V^{\threeD}:&=\{w\in H^1(D^{\threeD})|w_{|K}\in\mathbb{Q}_{N,d},\forall K\in\mathcal{T}^{\threeD}\},\\
V^{\threeD}_0:&=\{w\in V^{\threeD}|w_{|\partial D^{\threeD}_0\cup\partial D^{\threeD}_1}=0\},
\end{align}
\end{subequations}
where $\mathbb{Q}_{N,d}$ is the space of
$d$-variate polynomials of degree with respect to each variable at most $N$.
We use the $Q_N$ tensor-product Lagrange finite element
with the Gauss-Lobatto nodes in each direction.
We then consider the following problem whose solution is to converge to the
\emph{exact} solution in the physical domain given as:
\begin{align}
u^{\threeD}_{\textrm{exact}}=\frac{e^{-s\norm{x}}}{\norm{x}}:
\label{E:example_3d_ue}
\end{align}
\begin{subequations}
\begin{align}
\intertext{Find $u^{\threeD}\in V^{\threeD}$ such that}
\int_{D^{\threeD}}
\frac{1}{\prod\limits_{i\in\{0,1,2\}}g^{\threeD}_i}
\left(
s^2u^{\threeD}\ol{w^{\threeD}}
+\sum_{i\in\{0,1,2\}}(g^{\threeD}_i)^2u_{,i}^{\threeD}\ol{w_{,i}^{\threeD}}
\right)
\dif\bs{x}
&=0,\quad\forall w^{\threeD}\in V^{\threeD}_0,
\label{E:example_3d_domain}\\
{u^{\threeD}}_{|\partial D^{\threeD}_0}&=0,
\label{E:example_3d_end}\\
{u^{\threeD}}_{|\partial D^{\threeD}_1}&=\mathcal{I}(u^{\threeD}_{\textrm{exact}}),
\label{E:example_3d_dirichlet}
\end{align}
\label{E:example_3d}%
\end{subequations}
where, for $i\in\{1,2,3\}$:
\begin{align}
g^{\threeD}_i(x_i):=\left\{
\begin{array}{cl}
\gamma_0, &x_i\leq B_i,\\
\gamma_l, &B_i+(l-1)h<x_i\leq B_i+lh,\quad l\in\{1,\cdots,L\},
\end{array}
\right.
\end{align}
where $\gamma_0=1$ and:
\begin{align}
\gamma_l=\left(\cos\phi_l\!\cdot\!s + \frac{\sin^2\phi_l}{\cos\phi_l}\right)
\!\cdot\!\frac{1}{N + 1}\!\cdot\!h,
\quad \phi_l\in\left[0,\frac{\pi}{2}\right),
\quad l\in\{1,\cdots,L\}.
\label{E:example_3d_gamma_l}
\end{align}
In \eqref{E:example_3d_gamma_l}
the term in the parentheses was proposed in \cite{HagstromWarburton2009}
for the CRBCs for the case corresponding to $N=1$;
as optimisation of the ABC parameters is outside the scope of this paper,
we use this formula for all $N$,
simply choosing $\phi_l$ in $[0,\pi/2)$ according to the $L$-point
Gauss-Legendre quadrature.
The factors $1/(N+1)$ and $h$
were introduced specifically for the proposed ABC
to \emph{normalise} $\gamma_l$
across different values of $N$ and across different refinement levels,
respectively.
Although the standard symbol was used for the integral
in \eqref{E:example_3d_domain} for simplicity,
we only use the full integration rule in the physical domain, and
use the $N$-point reduced integration rule in $x_i>B_i$
in the $x_i$-direction for each $i\in\{1,2,3\}$.
\eqref{E:example_3d_end} is the termination condition for the ABC and
\eqref{E:example_3d_dirichlet} is the standard Dirichlet boundary condition,
where $\mathcal{I}:C^0(D^{\threeD})\rightarrow V^{\threeD}$ is the \emph{interpolation} operator from $C^0(D^{\threeD})$ to $V^{\threeD}$.
The homogeneous Neumann boundary condition
was weakly applied on $\partial D^{\threeD}_2$ in \eqref{E:example_3d_domain}.

We solved the above problem for $s\in\{+4.00+0.25\iu,+0.25+4.00\iu,+4.00\iu\}$.
For each $s$, we used four different finite element degrees,
$N\in\{1,2,3,4\}$, and
for each $N$, we considered four different refinement levels,
$\nref\in{1,2,3,4}$.
For each case, we varied the number of absorbing layers, $L$,
so that $LN\leq 12$.
We used MUMPS parallel sparse direct solver
\cite{mumps1,mumps2} via PETSc unless $(N,\nref)=(4,4)$.
For $(N,\nref)=(4,4)$, via PETSc,
we used the Generalized Minimal Residual method (KSPGMRES)
with the block Jacobi preconditioning (PCBJACOBI)
with the sub KSP type KSPPREONLY
and with the sub PC type
the incomplete factorization preconditioners (PCILU) or
the successive over relaxation (PCSOR).
The use of more elaborated preconditioners will be explored elsewhere.
For each case, we computed the relative $L^2$-norm error,
$\norm{u^{\threeD}-u^{\threeD}_{\textrm{exact}}}_2/\norm{u^{\threeD}_{\textrm{exact}}}_2$,
in the physical domain.
Fig.~\ref{Fi:example_3d_exact_0},
Fig.~\ref{Fi:example_3d_exact_1}, and
Fig.~\ref{Fi:example_3d_exact_2}
show the exact solutions for
$s=+4.00+0.25\iu$,
$s=+0.25+4.00\iu$, and
$s=+4.00\iu$.
Fig.~\ref{Fi:convergence_case_0},
Fig.~\ref{Fi:convergence_case_1}, and
Fig.~\ref{Fi:convergence_case_2},
show the plots of the computed relative $L^2$-norm errors for
$s=+4.00+0.25\iu$,
$s=+0.25+4.00\iu$, and
$s=+4.00\iu$, respectively,
for each $N$, for each $\nref$.
For each $s$, for each pair of $N$ and $\nref$,
we computed the interpolation error
$\norm{\mathcal{I}(u^{\threeD}_{\textrm{exact}})-u^{\threeD}_{\textrm{exact}}}_2/\norm{u^{\threeD}_{\textrm{exact}}}_2$
as an indicator of the discretisation error.
The interpolation error is shown with the black dotted line for each case.
The figures indicate that,
for each $s$, with $N$ and $\nref$ fixed,
the relative $L^2$-norm error decreases sharply to the discretisation error
as $L$ increases.
Such sharp convergence of the error is typically observed with PMDLs and CRBCs.

\section{Conclusion and future works}\label{S:conclusion}
We proposed a new class of ABCs for the Helmholtz equation,
which is parametrised by a tuple $(L,N)$.
The ABCs of type $(L,N)$
are derived from the layerwise constant PMLs with $L$ absorbing layers
using the $Q_N$ finite element
in conjunction with the $N$-point Gauss-Legendre quadrature reduced integration.
The proposed ABCs generalise the PMDLs \cite{GuddatiLim2006},
including them as type $(L,1)$,
which is known to outperform the standard PMLs and
to show equivalent performances to the CRBCs \cite{HagstromWarburton2009},
in terms of the accuracy and also in terms of the insensitivity
to the choice of the parameters.
The proposed ABCs of type $(L,N)$ inherit these favourable properties.
Like the PMDLs, the proposed ABCs only involve a single
Helmholtz-like equation over the entire computational domain,
even in existence of edges and corners.
As a generalisation of the PMDLs,
for any $N\in\mathbb{N}$,
the proposed ABCs allow for a monolithic discretisation of
the physical and the artificial domains
as long as a $Q_N$ finite element is used in the physical domain.
These features were validated in the one- and three-dimensional
numerical examples.
The analysis presented in this work were motivated by \cite{Ainsworth2004},
and appreciated the properties of the Jacobi polynomials.
We believe that this analysis gives more insight into the problems
and potentially advance the fields of ABCs in related areas.

We note that, though the ABCs of type $(L,N)$ and those of type $(LN,1)$
introduce the same number of auxiliary degrees of freedom
in the artificial domain and achieve the same order of accuracy,
the former only involve $L$ complex parameters while
the latter involve $LN$;
thus, one will have a finer control of the
absorption profile with the latter, which might be advantageous
if, e.g., one had a good \emph{a priori} knowledge of the target solution;
see also \cite{Guddati2016}.
Note also that
the matrix structure of the former is denser than that of the latter
in the artificial domain as the former use the $Q_N$ finite element and
the latter use the $Q_1$;
this could potentially be a drawback of using the former
if the size of the artificial domain was comparable with
that of the physical domain.

With regard to the proposed ABCs,
we have yet to investigate various aspects
that have been studied for
the PMDLs, the CRBCs, and the other related ABCs,
including applications to the wave equation \cite{GivoliNeta2003,HagstromWarburton2004,HagstromWarburton2009},
to the elastodynamics \cite{Baffet2012},
to the other elasticity problems \cite{SavadattiGuddati2012,SavadattiGuddati2012-1}, and 
to the problems with polygonal domains \cite{GuddatiLim2006}.
For the Helmholtz and the other equations,
we are also interested in using different families of finite elements
in conjunction with some reduced integration rules
in the artificial domain, so that
the ABCs can be applied with no added cost of implementation
when those families of elements are used for discretisation in the physical domain.
These ideas will be explored in the future works.

\section*{Acknowledgment}
This work was funded under the Engineering and Physical Sciences Research Council [grant numbers EP/R029423/1, EP/W029731/1].
I would like to thank Professor Ham for his consistent support to make this work possible.

\section*{Code availability}
The exact version of Firedrake used, along with all of the scripts employed in the experiments presented in this paper has been archived on Zenodo \cite{firedrake-zenodo}.

\bibliographystyle{amsplain}
\bibliography{reference}
% -- Appendix
%\appendix
%\section{Boundary conditions}\label{S:bc} 

\end{document}